\documentclass{article} 
\usepackage[margin=2.5cm]{geometry} 
\usepackage{a4wide}
\usepackage[T1]{fontenc} 
\usepackage{anyfontsize} 
\usepackage{ragged2e} 

\usepackage{amsfonts} 
\usepackage{amssymb, amsmath} 
\usepackage{amsthm} 
\usepackage{bbm} 
\numberwithin{equation}{section} 

\usepackage{comment}

\usepackage{graphicx} 
\usepackage{caption} 
\usepackage{subcaption} 
\usepackage{placeins} 
\usepackage{epstopdf} 
\usepackage{array,booktabs,multirow} 
\usepackage{booktabs} 
\usepackage{colortbl}
\usepackage{tikz}

\usetikzlibrary{arrows.meta,calc,positioning}

\usepackage{hyperref} 
\usepackage{xurl}      
\hypersetup{colorlinks=true, linkcolor=blue, filecolor=magenta, urlcolor=cyan, citecolor=teal,}
\newcommand{\citenote}[2]{\cite[#1]{#2}} 

\usepackage{xcolor} 
\definecolor{mygreen}{RGB}{28,172,0} 
\definecolor{mylilas}{RGB}{170,55,241} 
\usepackage{listings} 
\usepackage{extpfeil} 
\usepackage{cancel} 

\theoremstyle{definition}  
\newtheorem{theorem}{Theorem}[section]  
\newtheorem{remark}[theorem]{Remark}  
\newtheorem{proposition}[theorem]{Proposition}  
\newtheorem{corollary}[theorem]{Corollary}  

\begin{document}

\title{Data-driven calibration of a two-dimensional traffic flow model: Dimension reduction and optimization strategies}

\author{Simone G\"ottlich$^{1,*}$, \; Jacob Heieck$^{1}$}

\date{\today}

\maketitle

{\renewcommand{\thefootnote}{}%
\footnotetext{$^{1}$Department of Mathematics, University of Mannheim, 68131 Mannheim, Germany.}%
\footnotetext{$^{*}$\textbf{Corresponding author:} Simone G\"ottlich, Department of Mathematics, University of Mannheim, B6 28-29, 68159 Mannheim, Germany. Email: goettlich@uni-mannheim.de, Tel.: +49~621~181-2438}}

\maketitle

\begin{abstract}
We study the data-driven calibration of a two-dimensional first-order macroscopic traffic flow model, in which the movement along and across the road is described by a three-parameter flux model. The calibration is based on three real-world datasets, ranging from microscopic drone trajectories on German highways to aggregated urban sensor measurements. We recall the well-posedness of the underlying scalar conservation law and derive structural properties of the flux model. In particular, one parameter enters linearly and can be eliminated analytically, which reduces the non-linear least-squares problem from three to two dimensions. We solve the full and reduced calibration problem with two local optimization methods, Nelder-Mead and L-BFGS-B, and three global population-based methods Differential Evolution, Consensus-Based Optimization, and its deterministic variant Quasi-Monte Carlo consensus based optimization. In addition to accuracy and computational cost, we examine how sensitive the objective is to the parameters and how robust the algorithms are with respect to randomness and initialization. The reduced formulation is considerably cheaper and less sensitive to the choice of optimizer. The local methods are the most efficient when suitably initialized, whereas the global methods are more robust.
\end{abstract}



\section{Introduction}
\label{sec:introduction}

Macroscopic traffic flow models describe traffic in terms of aggregated quantities such as density, mean speed, and flux. Starting from the Lighthill-Whitham-Richards (LWR) model \cite{LWR,shock_waves_highway}, scalar conservation laws have become a standard tool for describing the formation and propagation of traffic waves \cite{traffic_flow_theory,traffic_flow_overview}. Their central constitutive ingredient is the fundamental diagram, i.e., the relation between density and flux. Instead of prescribing this relation a priori, it can be calibrated from trajectory or sensor data \cite{model_description_2,wang_macro_micro_calibration}. The growing availability of detailed traffic measurements makes this approach increasingly feasible.

On multi-lane roads, lane changes generate an additional traffic flux across the road. Macroscopic descriptions of this effect range from coupled lane-based balance laws to continuous two-dimensional models \cite{piu_herty_puppo_2023,balzotti_goettlich_2021}. In \cite{2D_traffic}, a two-dimensional first-order model was proposed and calibrated from trajectory data. There, the vehicle density $\rho=\rho(x,y,t)$ satisfies the scalar conservation law
\[
\partial_t\rho+\partial_x q^x(\rho)+\partial_y q^y(\rho)=0,
\]
where $x$ denotes the position along the road and $y$ the position across the lanes. The flux along the road belongs to the three-parameter family
\begin{align*}
    q_{\alpha,\lambda,p}(\rho)
=\alpha\Bigl((1-r)\,h_{\lambda,p}(0)+r\,h_{\lambda,p}(1)-h_{\lambda,p}(r)\Bigr),
\qquad
h_{\lambda,p}(r)=\sqrt{1+\lambda^2(r-p)^2},
\end{align*}
where $r=\rho/\rho_{\max}$ and $\rho_{\max}$ denotes the maximal density. The scale parameter $\alpha$ determines the magnitude of the flux, the critical-density parameter $p$ the density at which free-flow turns into congested traffic, and the curvature parameter $\lambda$ the sharpness of this transition.

We revisit this model using three datasets of substantially different character: drone trajectories from a highway section with an off-ramp and on-ramp \cite{ramp_dataset}, the highD trajectory dataset \cite{highD_dataset}, and aggregated one-dimensional sensor measurements of urban traffic in Heidelberg \cite{heidelberg_dataset}. This allows us to study the calibration across different road types, traffic regimes, and levels of aggregation.

The identification of a flux function from data is a non-linear inverse problem that has been studied from analytical and numerical viewpoints \cite{james_sepulveda_1999,li_evje_2023,li_geng_evje_2023}. In a parametric setting, a small fitting error does not guarantee that the parameters are uniquely determined. Different parameter combinations may yield almost identical fluxes, and an apparent failure of an optimizer may reflect a limitation of the model class rather than of the algorithm. We therefore combine a comparison of different optimization methods with an analysis of the parameterized flux family, which allows us to separate properties of the model from properties of the optimization methods. 

Our main contributions can be summarized as follows. We recall the well-posedness of the two-dimensional model in the Kružkov framework and establish structural properties of the underlying flux family. In particular, we show that the sign of $\lambda$ is not identifiable and that every non-degenerate flux is either concave or convex and has a constant sign. This latter observation motivates a modified flux in the transverse direction which, in contrast to \cite{2D_traffic}, allows for a change in the preferred lane-changing direction. Moreover, since $\alpha$ enters the flux linearly, it can be eliminated analytically, reducing the original three-parameter least-squares problem exactly to a two-dimensional optimization problem in $(\lambda,p)$. We further show that, for large $\lambda$, the flux depends essentially on $p$ and on the product $\alpha\lambda$, providing an explanation for the weak separate identifiability of $\alpha$ and $\lambda$. Finally, we compare the local methods Nelder-Mead and L-BFGS-B with the global population-based approaches Differential Evolution, Consensus-Based Optimization (CBO), and its deterministic variant QMC$^2$BO. The comparison is carried out on all three datasets for the full calibration problem and, in part, also for the reduced formulation. In addition, we investigate the sensitivity of the objective with respect to $(\lambda,p)$ and assess the robustness of the optimization methods with respect to randomness and initialization.

The remainder of the paper is organized as follows. Section~\ref{sec:data} introduces the datasets, their preprocessing, and the resulting fundamental diagrams. Section~\ref{sec:model} presents the traffic model and its analysis, including the dimension reduction of the calibration problem. Section~\ref{sec:optimization} describes the optimization algorithms, and Section~\ref{sec:numerical_experiments} reports the numerical results. Section~\ref{sec:conclusion} concludes the paper.

\section{Traffic datasets}
\label{sec:data}

This section introduces the three experimental traffic datasets used throughout this work. Although all datasets describe vehicular traffic in Germany, they differ with respect to the recording method and the level of aggregation.

The first two datasets, referred to as the \textit{On-/Off-ramp dataset} \cite{ramp_dataset}, published by the Institute of Highway Engineering (isac), and the \textit{highD dataset} \cite{highD_dataset}, published by the Institute for Automotive Engineering (ika), were recorded from drones and provide microscopic vehicle observations over several multi-lane highway segments. For every vehicle $i$ and time instant $t$, the recorded data include the position along the road ($x$-direction) and the position across the lanes ($y$-direction), i.e.,
$
(x_i(t),y_i(t)),
$
together with the corresponding velocity, i.e.,
$
(v_i^x(t),v_i^y(t)),
$
and additional vehicle-specific information such as accelerations and vehicle classes. Since both the motion along the road and across the lanes are available, these datasets are referred to as two-dimensional.

The third dataset, referred to as the \textit{Heidelberg sensor dataset} \cite{heidelberg_dataset}, was collected using a stationary traffic sensor operated by the City of Heidelberg as part of its urban traffic monitoring infrastructure. In contrast to the drone recordings, the sensors do not track individual vehicles along a road but record aggregated traffic at fixed locations. Consequently, only the motion along the road is available and no information across lanes can be extracted, resulting in a one-dimensional dataset.

In the following, we first describe the individual datasets and the quantities relevant for our numerical analysis. Afterwards, we explain how we transformed these different datasets into comparable macroscopic traffic data. Finally, we present a fundamental diagram for each dataset to provide an overview of their characteristics.

\subsection{Description of the datasets}
\label{sec:data_description}

\paragraph{\textbf{On-/Off-ramp dataset.}}

\begin{figure}[h]
    \centering
    \includegraphics[width=\linewidth]{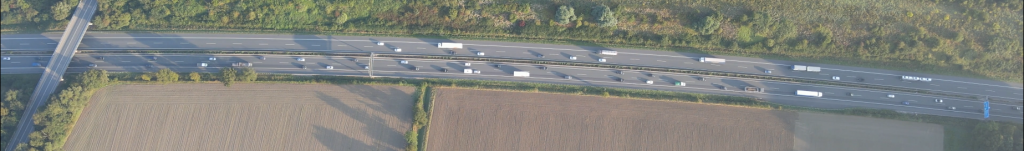}
    \includegraphics[width=\linewidth]{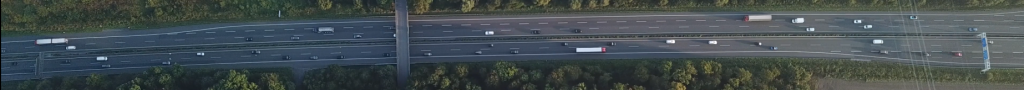}
    \caption{Aerial views of the highway sections contained in the On-/Off-ramp dataset \cite{ramp_dataset}. Top: western part. Bottom: eastern part.}
    \label{fig:road_ramp_dataset}
\end{figure}

The two-dimensional On-/Off-ramp dataset \cite{ramp_dataset} contains microscopic traffic observations obtained from aerial drone recordings of the German highway A43 near Münster. The highway consists of two lanes in each driving direction and includes a two-lane off-ramp in the direction from east to west and a one-lane on-ramp in the direction from west to east, see Figure~\ref{fig:road_ramp_dataset}. The two driving directions are provided as separate recordings.

For every observed vehicle $i$, the dataset contains the two-dimensional positions
$(x_i(t),y_i(t))$ and corresponding velocity $(v_i^x(t),v_i^y(t))$ observed every $0.1$ seconds. The recordings cover a $1200$ meters long highway section during the morning rush hour between 7:11 and 8:38 on September 6, 2021, resulting in a total observation period of $87$ minutes. Overall, the data contains trajectories of $5016$ vehicles on the off-ramp section and $3632$ vehicles on the on-ramp section. Although additional microscopic quantities are available, only the vehicle positions and velocities are required for the computation of density, flow and mean speed, which is conducted in the next subsection.

\paragraph{\textbf{highD dataset.}}

The highD dataset \cite{highD_dataset} is a large microscopic highway traffic dataset. Similar to the On-/Off-ramp dataset, it provides drone-based observations of the vehicle positions $(x_i(t),y_i(t))$ and corresponding velocity $(v_i^x(t),v_i^y(t))$ for every vehicle $i$.

\begin{figure}[h]
    \centering
    \includegraphics[width=\linewidth]{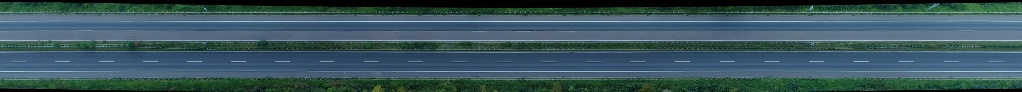}
    \caption{Aerial view of an example highway section contained in the highD dataset \cite{highD_dataset}.}
    \label{fig:road_highD_dataset}
\end{figure}

The recordings were conducted at six different highway locations in the Cologne area between 2017 and 2018. In total, the dataset consists of 60 recordings with an average duration of approximately $17$ minutes, covering road sections of roughly $420\,\mathrm{m}$ and highways with two to four lanes per driving direction, see for example Figure~\ref{fig:road_highD_dataset}. Both driving directions are recorded simultaneously within each recording. For our experiments, we use both directions of the first recording and pool them into a single dataset. The data is observed every $0.04$ seconds, giving a finer temporal resolution than the On-/Off-ramp dataset. Altogether, the dataset contains observations of more than $110\,000$ vehicles. As before, only the vehicle positions and velocity components will be used for the macroscopic traffic analysis.

\paragraph{\textbf{Heidelberg sensor dataset.}}

Unlike the previous two datasets, the Heidelberg sensor dataset \cite{heidelberg_dataset} consists of sensor-based traffic data. Instead of tracking individual vehicles over a road section, stationary sensors installed at 50 individual lanes throughout the urban road network in Heidelberg continuously record the number of passing vehicles together with their average speed.

The data has been collected since December 18, 2024, and are aggregated observations over consecutive $15$-minute intervals. Since each sensor only monitors a single traffic lane, no vehicle trajectories or changing lane information are available. 

\begin{figure}[h]
    \centering
    \includegraphics[width=0.8\linewidth]{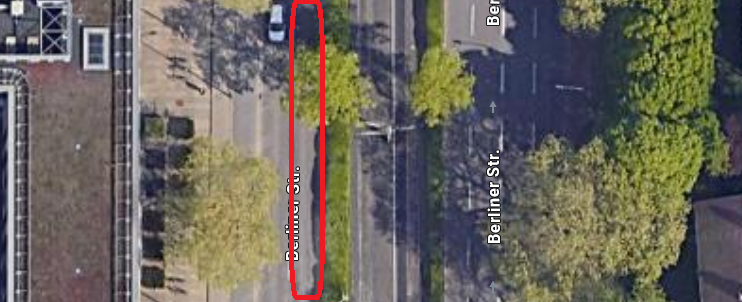}
    \caption{Aerial view of the selected lane at \emph{Berliner Straße / Im Neuenheimer Feld} in the Heidelberg sensor dataset \cite{heidelberg_dataset}.}
    \label{fig:road_heidelberg_dataset}
\end{figure}

For the numerical experiments, we select the sensor located at \emph{Berliner Straße / Im Neuenheimer Feld} between January and April 2025, see Figure~\ref{fig:road_heidelberg_dataset}. During this period, the sensor recorded approximately $2\,260\,000$ vehicles, making it one of the busiest sensor locations in the available network and therefore particularly suitable for calibrating macroscopic traffic models. In contrast to the 2-dimensional drone datasets, the sensor data already provides aggregated traffic counts $q(t)$ and average vehicle speeds $u^x(t)$, requiring a different preprocessing procedure that will be described in the next subsection.

\medskip

To summarize the main characteristics of all three datasets, we present Table~\ref{tab:datasets}. For the highD dataset, we consider the representative first recording, which is also used throughout the numerical experiments presented in this work.

\begin{table}[ht]
\centering
\caption{Comparison of the three traffic datasets considered in this work.}
\small
\renewcommand{\arraystretch}{1.7}
\begin{tabular}{p{3.0cm}p{3.7cm}p{3.6cm}p{3.8cm}}
\toprule
\textbf{Property}
&
\textbf{On-/Off-ramp data}
&
\textbf{highD data}
&
\textbf{Heidelberg sensor data}
\\
\midrule

Method
&
Drone-based (microscopic)
&
Drone-based (microscopic)
&
Fixed sensor (macroscopic)
\\

Location
&
Highway A43 near Münster
&
Highway near Cologne
&
Urban road in Heidelberg
\\

Road length
&
$1200\,\mathrm{m}$
&
$430.2\,\mathrm{m}$
&
Fixed sensor location
\\

Recording period
&
09/06/2021 07:11-08:38
&
09/2017 08:38-08:53
&
01/2025 - 04/2025
\\

Driving directions
&
Separated into 2 recordings
&
Compiled in 1 recording
&
Single driving direction
\\

Number of lanes
&
2 lanes per direction
&
2 lanes per direction
&
1 lane for given direction
\\

Sampling interval
&
$0.1\,\mathrm{s}$
&
$0.04\,\mathrm{s}$
&
$15\,\mathrm{min}$
\\

Observed data
&
$x_i(t)$, $y_i(t)$, $v_i^x(t)$, $v_i^y(t)$
&
$x_i(t)$, $y_i(t)$, $v_i^x(t)$, $v_i^y(t)$
&
$q^x(t)$, $u^x(t)$
\\

\bottomrule
\end{tabular}

\label{tab:datasets}
\end{table}

\subsection{Data preparation}
\label{sec:data_preparation}

The previous subsection introduced the three traffic datasets considered in this work. We now describe how the recorded observations are transformed into the macroscopic traffic data required for fundamental diagrams. In particular, we compute the density $\rho$, the traffic flux $q$, and the mean speed $u$, measured in vehicles per kilometer, vehicles per hour, and kilometers per hour, respectively.

Since the On-/Off-ramp and highD datasets provide microscopic vehicle observations, while the Heidelberg sensor data already consist of aggregated traffic observations, the preprocessing differs between these two types of data. We therefore first describe the general procedure for the microscopic datasets before discussing the sensor-based data separately. The preprocessing of the microscopic data closely follows the procedure proposed in \cite{2D_traffic}. For a more detailed discussion of the conversion from microscopic vehicle trajectories to macroscopic traffic quantities, we refer to \cite{2D_traffic,traffic_flow_theory,traffic_flow_overview}.

In the case of microscopic models, we consider two dimensions, since both the motion along the road and the movement across lanes are observed. Consequently, the described procedure is applied separately to the $x$- and $y$-directions. For the sensor-based dataset, only the traffic flow along the road is available, so no preprocessing of the $y$-direction is required.

\paragraph{\textbf{Microscopic highway data.}}

For the On-/Off-ramp and highD datasets, the available observations consist of the two-dimensional vehicle positions $(x_i(t),y_i(t))$ and corresponding velocities $(v_i^x(t),v_i^y(t))$ for every vehicle $i$ and time $t$. Since the two datasets were recorded with different temporal resolutions ($0.1\,\mathrm{s}$ and $0.04\,\mathrm{s}$), the first preprocessing step consists of transforming all observations onto a common one-second temporal grid. More precisely, every observation is assigned to
$
t_k=\lfloor t\rfloor.
$

For every discrete time instant $t_k$, let $N(t_k)$ denote the number of vehicles within the observed road section of length $L$. The density is then given by
\begin{equation}
\tilde{\rho}(t_k) = \frac{N(t_k)}{L}, \qquad k=0,\ldots,M,
\end{equation}
where $L$ is measured in kilometers and $M$ denoted the total number of time instants.

The mean velocity is obtained by averaging the microscopic vehicle velocities in each direction
\begin{equation}
\tilde{u}^x(t_k) = \frac{1}{N(t_k)} \sum_{i=1}^{N(t_k)} v_i^x(t_k), \qquad
\tilde{u}^y(t_k) = \frac{1}{N(t_k)} \sum_{i=1}^{N(t_k)} v_i^y(t_k).
\end{equation}
The corresponding traffic fluxes follow directly from the hydrodynamic relation
\begin{equation}
\tilde{q}^x(t_k) = \tilde{\rho}(t_k)\tilde{u}^x(t_k), \qquad 
\tilde{q}^y(t_k) = \tilde{\rho}(t_k)\tilde{u}^y(t_k).
\end{equation}

The one-second time series captures microscopic traffic dynamics, but still shows short-term fluctuations. To obtain smoother macroscopic quantities suitable for the computation of fundamental diagram, we apply a moving average over $m$ consecutive one-second observations. For the density and average fluxes, this yields
\begin{equation}\label{eq:aggregation_1}
\rho_{k_0} = \frac{1}{m} \sum_{k=k_0}^{k_0+m-1} \tilde{\rho}(t_k),\qquad
q_{k_0}^x = \frac{1}{m} \sum_{k=k_0}^{k_0+m-1} \tilde{q}^x(t_k), \qquad
q_{k_0}^y = \frac{1}{m} \sum_{k=k_0}^{k_0+m-1} \tilde{q}^y(t_k)  
\end{equation}
for $k_0=0,\ldots,M+1-m$. Finally, the averaged mean velocities are recovered via
\begin{equation}\label{eq:aggregation_2}
u_{k_0}^x = \frac{q_{k_0}^x}{\rho_{k_0}}, \qquad
u_{k_0}^y = \frac{q_{k_0}^y}{\rho_{k_0}}.
\end{equation}

Throughout this work, we choose
$   
m=180,
$  
corresponding to a three-minute aggregation interval. The resulting data pairs $(\rho_{k_0},q^x_{k_0})$, $(\rho_{k_0},q^y_{k_0})$, $(\rho_{k_0},u^x_{k_0})$ and $(\rho_{k_0},u^y_{k_0})$ form the fundamental diagrams and speed-density diagrams for each direction, which are discussed in the next subsection.

\paragraph{\textbf{Sensor-based urban data.}}

The Heidelberg sensor data require a different preprocessing procedure. Each sensor records the total number of passing vehicles $q^x(t)$ together with their average speed $u^x(t)$ over consecutive intervals of $15$ minutes. The corresponding traffic flux is directly obtained by scaling these traffic counts into an hourly flow rate, i.e.
\begin{equation}
    \tilde{q}^x(t_k)=4q^x(t_k),
\end{equation}
while the mean speed is directly provided by the sensor, i.e.
\begin{equation}
   \tilde{u}^x(t_k)=u^x(t_k).
\end{equation}
The density is then computed from the hydrodynamic relation
\begin{equation}
    \tilde{\rho}(t_k)=\frac{\tilde{q}^x(t_k)}{\tilde{u}^x(t_k)}.
\end{equation}

Since the sensor data are already aggregated over comparatively long time intervals, an additional moving average over only 
$ m=4$  
consecutive time instants $t_k$ identical to \eqref{eq:aggregation_1} and \eqref{eq:aggregation_2} is applied. This corresponds to a one-hour aggregation interval. Again, the resulting data pairs $(\rho_{k_0},q^x_{k_0})$ and $(\rho_{k_0},u^x_{k_0})$ form the fundamental diagrams and speed-density diagrams for the $x$-direction. The discussion is done in the next subsection

\subsection{Data visualization}
\label{sec:data_visualization}

Based on the macroscopic quantities obtained in the previous subsection, we now construct the corresponding \emph{fundamental diagrams} and \emph{speed-density diagrams}. More precisely, the data pairs $(\rho_{k_0},q^x_{k_0})$, $(\rho_{k_0},q^y_{k_0})$, $(\rho_{k_0},u^x_{k_0})$ and $(\rho_{k_0},u^y_{k_0})$ are used to visualize the dependence of the traffic flux and the mean speed on the vehicle density. Throughout this work, the density is shown on the horizontal axis, while the vertical axis displays either the traffic flux or the mean speed.

These diagrams constitute one of the fundamental tools in macroscopic traffic flow theory. They provide valuable insight into the transition between free-flow traffic, congested traffic, and, in extreme cases, bumper-to-bumper traffic. For a detailed discussion of these diagrams, we refer to, e.g., \cite{2D_traffic,fundamendal_diagrams_physics_of_traffic,fundamental_diagrams_best,fundamental_diagrams_75_years}.

\medskip

\begin{figure}[ht]
    \centering
    \includegraphics[width=0.85\textwidth]{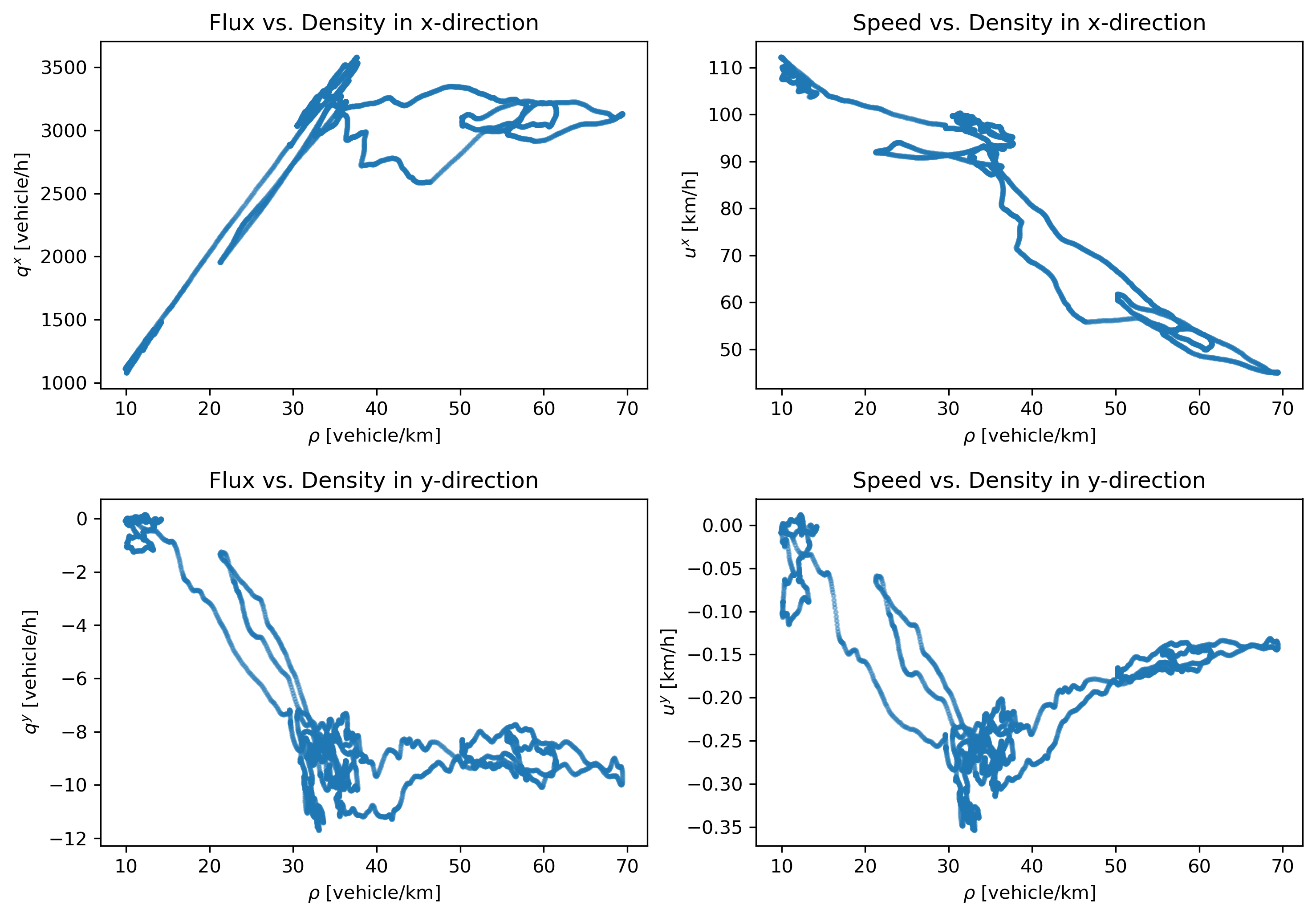}
    \caption{Fundamental diagrams (left) and speed-density diagrams (right) obtained from the Off-ramp dataset. The upper row shows the movement along the road ($x$-direction), while the lower row illustrates the movement across lanes ($y$-direction).}
    \label{fig:fundamental_diagrams_offramp}
\end{figure}

We first consider the Off-ramp dataset. Figure~\ref{fig:fundamental_diagrams_offramp} shows the corresponding fundamental diagrams (left) and speed-density diagrams (right). The upper row illustrates the movement along the road ($x$-direction), while the lower row displays the movement across lanes ($y$-direction).

The diagrams in $x$-direction show the classical behavior observed in free-flow traffic. Up to approximately $40$ vehicles per kilometer, the flux increases almost linearly with the density, indicating free-flow. Since the data were recorded during the morning rush hour, the data also contain several congested traffic states, where the density continues to increase while the traffic flux is stable or decreases. The maximum observed density is approximately $70$ vehicles per kilometer and bumper-to-bumper traffic is not observed.

The diagrams describing the movement in $y$-direction show a different behavior. The flux and the mean speed are much smaller than their counterparts in $x$-direction. This is because velocities caused by changing lanes are expected to be much smaller than velocities caused by movement along the road. Nevertheless, these quantities should not be neglected. As demonstrated in \cite{2D_traffic}, changing lanes has a measurable influence on the overall traffic evolution, since approximately $15\%$ of all vehicles change lanes while being on the German highway section they analyzed.

A second difference of the movement in $y$-direction is the sign. While speeds and fluxes in $x$-direction are strictly non-negative, since vehicles travel only in the driving direction, the quantities for the movement across lanes can attain positive and negative values. Throughout this work, positive values correspond to lane changes towards the leftmost lane, whereas negative values indicate lane changes towards the rightmost lane. Since the considered road section contains an off-ramp, most lane changes are directed towards the right, resulting in mostly negative speeds and fluxes.

\medskip

\begin{figure}[h]
    \centering
    \includegraphics[width=0.85\textwidth]{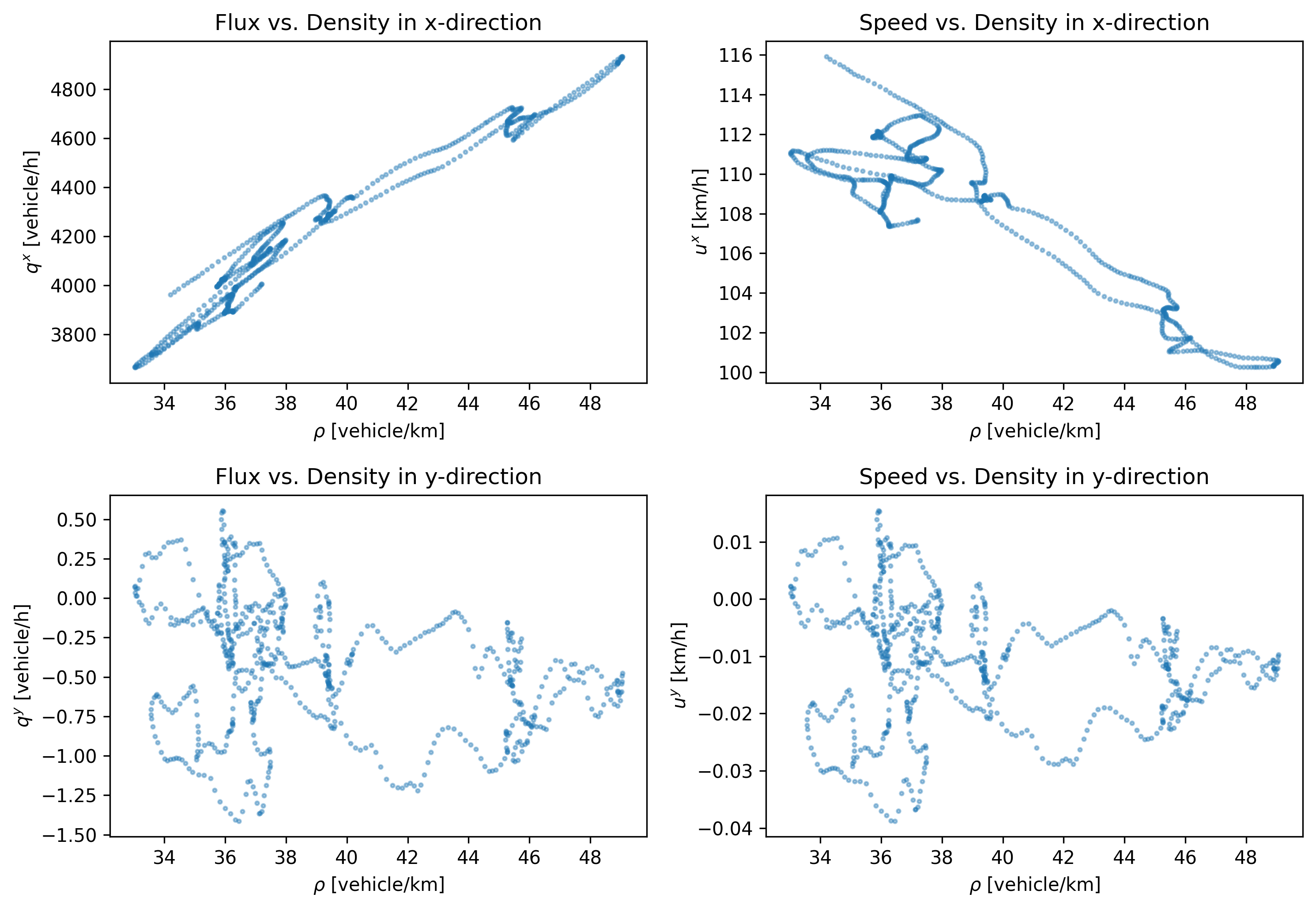}
    \caption{Fundamental diagrams (left) and speed-density diagrams (right) for the highD dataset for the $x$-direction (movement along the road) (top) and the $y$-direction (movement across lanes) (bottom).}
    \label{fig:fundamental_diagrams_highD}
\end{figure}

Next, we consider the first recording of the highD dataset shown in Figure~\ref{fig:road_highD_dataset}. In contrast to the On-/Off-ramp dataset, both driving directions are contained in the same recording. The fundamental diagrams are shown in Figure \ref{fig:fundamental_diagrams_highD}.

Compared with the Off-ramp dataset, the diagrams contain fewer data points. The diagrams in $x$-direction again show the classical behavior of free-flow traffic. In contrast to the Off-ramp data, this highway section remained uncongested throughout the recording.

The diagrams in $y$-direction contain both positive and negative values, indicating lane changes in both directions. A slightly larger proportion of negative values suggests that lane changes towards the rightmost lane occur more frequently. This behavior is consistent with the observed free-flow traffic, where vehicles typically return to the right lane after completing overtaking maneuvers.

\begin{figure}[ht]
    \centering
    \includegraphics[width=\textwidth]{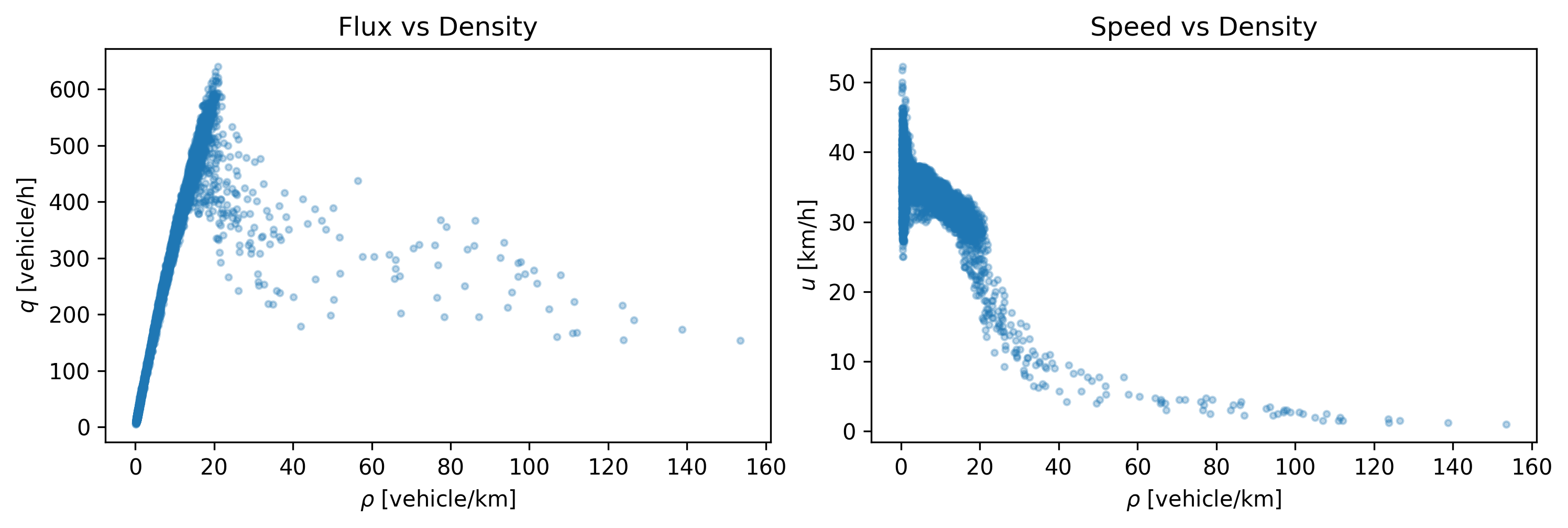}
    \caption{Fundamental diagram (left) and speed-density diagram (right) for the Heidelberg sensor dataset.}
    \label{fig:fundamental_diagrams_heidelberg}
\end{figure}

Finally, Figure~\ref{fig:fundamental_diagrams_heidelberg} shows the corresponding diagrams for the Heidelberg sensor dataset in $x$-direction, since only traffic quantities in this direction are available.

The resulting diagrams differ from those of the microscopic highway datasets. Due to the longer observation period of approximately four months, the diagrams contain much larger and smoother observations, resulting in macroscopic trends with only a few isolated outliers.

The fundamental diagram again shows an approximately monotone increase of the traffic flux up to densities of about $20$ vehicles per kilometer. Beyond this point, increasing densities yield decreasing fluxes, indicating congested traffic. The speed-density diagram provides a clear illustration of this transition. While vehicles travel close to the speed limit of $30$ km/h under $20$ vehicles per kilometer, the mean speed continuously decreases as the density increases beyond this value. At densities above approximately $80$ vehicles per kilometer, the speeds approach zero, indicating periods of bumper-to-bumper traffic. Intermediate densities between approximately $30$ and $80$ vehicles per kilometer correspond to moderately congested traffic, where vehicles continue to move but at reduced speeds. We note that only two data points exceed the maximal density $\rho_{\max}\approx133$ vehicles per kilometer, which indicates that $\rho_{\max}$ provides a reasonable estimate for this dataset. We attribute these two observations to the fact that $\rho=\tilde{q}^x/\tilde{u}^x$ is numerically sensitive to measurement noise whenever the measured speed is close to zero.

In general, despite the different observation periods, road characteristics, and aggregation levels, all three datasets show similar behavior. The main differences arise from the presence of highway ramps, urban traffic, or lane-changing dynamics, rather than from the underlying traffic behavior itself.

\begin{remark}
    These diagrams are obtained using the aggregation intervals introduced in Section~\ref{sec:data_preparation}, namely a three-minute moving average for the microscopic highway data and a one-hour moving average for the Heidelberg sensor data. Changing the aggregation window mainly affects the spread of the data points. The overall structure of the diagrams remains essentially unchanged. Additional examples for different aggregation windows can be found in the accompanying GitHub repository \cite{git_2D_data}.
\end{remark}

\section{Macroscopic traffic model}
\label{sec:model}

We now replace these discrete observations by a continuous mathematical model that captures the underlying traffic dynamics using only a small number of parameters. Once such a model has been chosen, we formulate the corresponding parameter identification problem.

Throughout this work, we employ the two-dimensional first-order macroscopic traffic model proposed in \cite{2D_traffic}. Rather than introducing a new flux family, we use this model as a common benchmark and study its structural properties and calibration from real traffic data. Related macroscopic descriptions of lane-changing and multi-lane effects include the two-dimensional multi-class model \cite{balzotti_goettlich_2021} and coupled multi-lane balance laws derived from microscopic dynamics \cite{piu_herty_puppo_2023}. For a detailed motivation of the particular model considered here, its relation to the classical one-dimensional Lighthill-Whitham-Richards (LWR) model, and its experimental validation, we refer to \cite{2D_traffic,LWR,shock_waves_highway}.

\subsection{Two-dimensional macroscopic traffic model}
\label{sec:2d_model}

Following \cite{2D_traffic}, we describe the evolution of the vehicle density by a two-dimensional scalar conservation law. Let
\[
\rho=\rho(x,y,t)
\]
denote the vehicle density at position $(x,y)$ and time $t$, where the $x$-direction corresponds to the movement along the road and the $y$-direction describes the movement across the lanes. The density satisfies
\begin{equation}
\partial_t \rho + \nabla\cdot F(\rho)=0,
\qquad t>0,\quad (x,y)\in\Omega,
\label{eq:lwr_2d}
\end{equation}
where 
\[
\Omega=[0,L^x]\times[0,L^y]
\]
denotes the considered road section and
\[
F(\rho)=
\begin{pmatrix}
q^x(\rho)\\
q^y(\rho)
\end{pmatrix}
\]
is the two-dimensional traffic flux. Equivalently, \eqref{eq:lwr_2d} reads
\[
\partial_t\rho
+\partial_x q^x(\rho)
+\partial_y q^y(\rho)=0.
\]

The two components of the flux describe different mechanisms. The flux $q^x$ accounts for the movement of vehicles along the road, whereas the flux $q^y$ represents their movement across the road due to lane changes. Both are related to the corresponding mean velocities through the fundamental relations
\begin{equation}
q^x(\rho)=\rho\,u^x(\rho),
\qquad
q^y(\rho)=\rho\,u^y(\rho).
\label{eq:hydrodynamic_relation}
\end{equation}

The velocity along the road $u^x$ is assumed to decrease as the vehicle density increases. In particular,
\[
u^x(\rho_{\max})=0,
\]
where $\rho_{\max}>0$ denotes the maximum vehicle density corresponding to bumper-to-bumper traffic. Consequently, the flux $q^x$ vanishes in the absence of vehicles and at maximum density.

The behavior in $y$-direction is considerably more involved, since lane-changing decisions depend on traffic density, driver preferences, and traffic regulations. In free-flow situations, vehicles tend to move to the rightmost lane. However, as traffic density increases, lane changes become increasingly difficult. Consequently, unlike the behavior in $x$-direction, the velocity in $y$-direction may attain positive and negative values.

We emphasize that \eqref{eq:lwr_2d} describes traffic as a continuous two-dimensional flow. Although the model captures lane-changing behavior, it should not be interpreted as a discrete multi-lane model, where individual lanes are represented separately. Instead, the movement across lanes is treated as a continuous variable.

Before turning to the parameterized flow model, we briefly discuss the well-posedness of the underlying traffic model. Since the experimental road section $\Omega$ is bounded, a complete analytical treatment on $\Omega$ would require the specification of appropriate entropy boundary conditions. As these boundary conditions are not part of the parameter identification considered below, we formulate the standard Cauchy problem on $\mathbb{R}^2$ \cite{kruzkov1970}. Note that the calibration considered in this work only involves the flux functions $q^x$ and $q^y$ and does not require solving \eqref{eq:lwr_2d} on $\Omega$. Hence, no boundary conditions enter the numerical experiments.

\begin{proposition}[Well-posedness of the two-dimensional traffic model]
\label{prop:wellposedness}
Let $F\in C^1([0,\rho_{\max}];\mathbb{R}^2)$ and let
\[
\rho_0\in L^1(\mathbb{R}^2)\cap L^\infty(\mathbb{R}^2),
\qquad
0\leq \rho_0\leq \rho_{\max}\quad\text{a.e. in }\mathbb{R}^2.
\]
Then the Cauchy problem
\begin{equation}
\label{eq:cauchy_traffic}
\partial_t\rho+\nabla\cdot F(\rho)=0,
\qquad
\rho(\cdot,0)=\rho_0,
\end{equation}
admits a unique Kružkov entropy solution $\rho\in L^\infty(\mathbb{R}^2\times(0,\infty))$. Moreover,
\[
0\leq \rho(x,y,t)\leq \rho_{\max}
\quad\text{for a.e. }(x,y)\in\mathbb{R}^2,\ t>0.
\]
If $\rho$ and $\widetilde\rho$ are the entropy solutions corresponding to initial data $\rho_0$ and $\widetilde\rho_0$, respectively, then
\begin{equation}
\label{eq:l1_contraction}
\|\rho(\cdot,t)-\widetilde\rho(\cdot,t)\|_{L^1(\mathbb{R}^2)}
\leq
\|\rho_0-\widetilde\rho_0\|_{L^1(\mathbb{R}^2)}
\qquad\text{for all }t\geq0.
\end{equation}
\end{proposition}

\begin{proof}
Since $F\in C^1([0,\rho_{\max}];\mathbb{R}^2)$, it can be extended to a function $\widetilde F\in C^1(\mathbb{R};\mathbb{R}^2)$ with bounded derivative (e.g., by a $C^1$ extension that is affine outside a neighborhood of $[0,\rho_{\max}]$). The classical Kružkov theory for scalar conservation laws in several space dimensions \cite{kruzkov1970} therefore yields existence and uniqueness of an entropy solution to \eqref{eq:cauchy_traffic}, together with the $L^1$-contraction property \eqref{eq:l1_contraction} \citenote{\text{Theorem 1}}{kruzkov1970}. The constant functions $\rho\equiv0$ and $\rho\equiv\rho_{\max}$ are entropy solutions. The comparison principle \citenote{\text{Theorem 3}}{kruzkov1970}, applied to the initial bounds $0\leq\rho_0\leq\rho_{\max}$, consequently gives
\[
0\leq\rho(\cdot,t)\leq\rho_{\max}
\]
for all $t\geq0$. Hence the solution remains in $[0,\rho_{\max}]$ and is independent of the particular extension of $F$ outside $[0,\rho_{\max}]$.
\end{proof}

Proposition~\ref{prop:wellposedness} shows that, once the constitutive relation $F$ is prescribed, the macroscopic traffic dynamics are well defined and stable with respect to perturbations of the initial density. For the data-driven setting of this work, however, the flux itself is inferred from measurements. Our main focus is therefore on identifying the constitutive relation
\[
\rho\longmapsto F(\rho)
\]
from experimental traffic data. We next introduce the parameterized family of flux functions used throughout the calibration procedure.

\subsection{Parameterized flow model}
\label{sec:data_fitting}

Following \cite{2D_traffic}, we consider a three-parameter family of flux functions. We first introduce this family for the $x$-direction.

Let
\begin{equation}
q_{\alpha,\lambda,p}(\rho) =\alpha\left(d_1+(d_2-d_1)\frac{\rho}{\rho_{\max}} -\sqrt{1+d_3(\rho)^2}\right),
\label{eq:flux_param}
\end{equation}
where
\begin{equation*}
d_1=\sqrt{1+(\lambda p)^2},\qquad
d_2=\sqrt{1+(\lambda(1-p))^2},\qquad
d_3(\rho)=\lambda\left(\frac{\rho}{\rho_{\max}}-p\right),
\end{equation*}
and set 
\begin{equation}
    q^x(\rho) :=q_{\alpha^x,\lambda^x,p^x}(\rho).
\end{equation} 
The maximal density $\rho_{\max}$ is determined by the ratio between the number of lanes and the typical vehicle length of 5 meters, plus 50\% of additional safety distance, so that
\begin{equation*}
\rho_{\max} = \frac{1000 \text{ m}\times \text{\# lanes}}{7.5\text{ m}}.
\end{equation*}
This corresponds to the theoretical bumper-to-bumper density of the considered road. For the $x$-direction, \eqref{eq:flux_param} coincides with the model proposed in \cite{2D_traffic}.

Details concerning the derivation, motivation, and validation of the model \eqref{eq:flux_param} can be found in \cite{model_description_1,model_description_2,2D_traffic}. In this work, we focus on the interpretation of the three parameters $\alpha^x$, $p^x$, and $\lambda^x$ in terms of the resulting fundamental diagram. We first give a brief, informal overview of how each parameter shapes the fitted diagram, and then make this interpretation mathematically precise in Section \ref{sec:parameter_identification}.

\textbf{Capacity-scaling parameter $\alpha^x$.}
The parameter $\alpha^x\in\mathbb{R}$ acts as a multiplicative scaling factor for the traffic flux. For fixed values of $p^x$ and $\lambda^x$, increasing $\alpha^x$ increases the flux throughout the fundamental diagram and, in particular, its maximum value. Therefore, roads with high traffic capacity, such as multi-lane highways, are expected to have larger $\alpha^x$ than roads with lower capacities, for example urban roads.

\textbf{Critical-density parameter $p^x$.}
The parameter $p^x$ mainly controls the density at which the transition from free-flow to congested traffic occurs. It describes how much traffic density can be achieved before drivers begin to reduce their speed significantly. As $p^x$ is intended to denote a dimensionless reference density relative to the maximal density $\rho_{\max}$, we restrict it to
\begin{equation*}
p^x\in[0,1].
\end{equation*}

\textbf{Curvature parameter $\lambda^x$.}
The dimensionless parameter $\lambda^x$ controls the curvature around $p^x$. It states how abruptly the transition from free-flow to congestion occurs. Large values of $\lambda^x$ indicate a sharper transition of the flux function. 

In addition to this interpretation, the following proposition summarizes useful properties of this parameterized flux family.

\begin{proposition}[Properties of the parameterized flow model]
\label{prop:flux_structure}
Let $\rho_{\max}>0$ and let $q_{\alpha,\lambda,p}$ be defined by \eqref{eq:flux_param}. Then the following properties hold.
\begin{enumerate}
    \item The flux satisfies
    \begin{equation}
    q_{\alpha,\lambda,p}(0)=q_{\alpha,\lambda,p}(\rho_{\max})=0.
    \label{eq:flux_boundary_values}
    \end{equation}
    
    \item The parameterization is invariant under a change of sign of $\lambda$, i.e.,
    \begin{equation}
    q_{\alpha,\lambda,p}(\rho)=q_{\alpha,-\lambda,p}(\rho)
    \qquad\text{for all }\rho\in[0,\rho_{\max}].
    \label{eq:lambda_symmetry}
    \end{equation}
    Consequently, the sign of $\lambda$ cannot be identified from data.
    
    \item If $\lambda\neq0$, then
    \begin{equation}
    \frac{\partial^2 q_{\alpha,\lambda,p}(\rho)}{\partial\rho^2}
    = -\frac{\alpha\lambda^2}{\rho_{\max}^2
    \left(1+\lambda^2\left(\frac{\rho}{\rho_{\max}}-p\right)^2\right)^{3/2}}.
    \label{eq:flux_second_derivative}
    \end{equation}
    Hence $q_{\alpha,\lambda,p}$ is strictly concave for $\alpha>0$ and strictly convex for $\alpha<0$.
    
    \item If $\lambda\neq0$, the sign of the flux is determined by the sign of $\alpha\neq0$
    \begin{equation}
    \operatorname{sign}\bigl(q_{\alpha,\lambda,p}(\rho)\bigr)
    =\operatorname{sign}(\alpha),
    \qquad 0<\rho<\rho_{\max}.
    \label{eq:constant_sign}
    \end{equation}
\end{enumerate}
\end{proposition}

\begin{proof}
Set $r=\frac{\rho}{\rho_{\max}}$ and define
\[
h_{\lambda,p}(r)=\sqrt{1+\lambda^2(r-p)^2},\qquad r\in[0,1].
\]
Then $d_1=h_{\lambda,p}(0)$, $d_2=h_{\lambda,p}(1)$, and
\begin{equation}
q_{\alpha,\lambda,p}(\rho)
=\alpha\bigl((1-r)h_{\lambda,p}(0)+r h_{\lambda,p}(1)-h_{\lambda,p}(r)\bigr).
\label{eq:flux_chord_representation}
\end{equation}
The boundary identities in \eqref{eq:flux_boundary_values} follow immediately by setting $r=0$ and $r=1$. Moreover, $h_{\lambda,p}$ depends on $\lambda$ only through $\lambda^2$, which proves \eqref{eq:lambda_symmetry}.

A direct calculation gives
$$
h_{\lambda,p}'(r) = \frac{\lambda^2(r-p)} {\sqrt{1+\lambda^2(r-p)^2}}, \qquad h_{\lambda,p}''(r) = \frac{\lambda^2} {\left(1+\lambda^2(r-p)^2\right)^{3/2}}.
$$
Hence, $h_{\lambda,p}$ is strictly convex for $\lambda\neq0$.

Differentiating \eqref{eq:flux_chord_representation}, using $r=\frac{\rho}{\rho_{\max}}$ and the chain rule yields
\begin{align*}
\frac{\partial q_{\alpha,\lambda,p}(\rho)}{\partial\rho}
&= \alpha\left( \frac{1}{\rho_{\max}} \left(h_{\lambda,p}(1)-h_{\lambda,p}(0)\right) - \frac{1}{\rho_{\max}} h_{\lambda,p}'(r) \right).
\end{align*}
Differentiating once more gives
\begin{align*}
\frac{\partial^2 q_{\alpha,\lambda,p}(\rho)}{\partial\rho^2} = -\frac{\alpha}{\rho_{\max}^2}
h_{\lambda,p}''(r).
\end{align*}
Substituting the expression for $h_{\lambda,p}''$, we obtain \eqref{eq:flux_second_derivative}.

Since $\rho_{\max}>0$ and $\lambda\neq0$, the denominator is strictly positive. Consequently, the sign of the second derivative is determined by $-\alpha$. Thus, $q_{\alpha,\lambda,p}$ is strictly concave for $\alpha>0$ and strictly convex for $\alpha<0$.

Finally, strict convexity of $h_{\lambda,p}$ implies that, for $0<r<1$,
\[
h_{\lambda,p}(r)
< r h_{\lambda,p}(1) + (1-r)h_{\lambda,p}(0).
\]
Together with \eqref{eq:flux_chord_representation}, this proves \eqref{eq:constant_sign}.
\end{proof}

Proposition \ref{prop:flux_structure} has two immediate consequences for the following calibration problem. First, the symmetry \eqref{eq:lambda_symmetry} implies that allowing positive and negative values of $\lambda$ leads to non-unique solutions of the problem. In particular, if a value $\lambda^\star$ is optimal for the problem, then $-\lambda^\star$ yields exactly the same flux and is therefore optimal as well. To avoid this non-uniqueness, we may restrict the curvature parameter to $\lambda\geq0$ without loss of generality.

Second, the sign property \eqref{eq:constant_sign} reveals a limitation of this modeling choice for the $y$-direction. For every fixed parameter vector $(\lambda,\alpha,p)$ with $\lambda\neq0$, the flux $q_{\alpha,\lambda,p}$ has a constant sign on $(0,\rho_{\max})$. For the $x$-direction, this is a natural assumption. We expect all vehicles to travel in the same direction, and therefore $\alpha^x\geq0$ directly gives a reasonable $x$-flux.

For the $y$-direction, movement across lanes may occur in both directions depending on the traffic situation. A model with constant sign on $(0,\rho_{\max})$ cannot capture this, and \eqref{eq:flux_param} is therefore not directly suitable for $y$. We instead multiply $q_{\alpha,\lambda,p}$ by the signed factor $\frac{\rho}{\rho_{\max}}-p$, which, for $p\in[0,1]$, introduces at most one sign change, located at $\rho=p\,\rho_{\max}\in[0,\rho_{\max}]$, without otherwise altering the qualitative shape of the diagram. This motivates the $y$-flux model
\begin{equation}
q^y_{\alpha,\lambda,p}(\rho)
:=\left(\frac{\rho}{\rho_{\max}}-p\right)q_{\alpha,\lambda,p}(\rho),
\label{eq:flux_param_y}
\end{equation}
with its own parameter set $(\alpha^y,\lambda^y,p^y)$, where, as for $p^x$, we restrict $p^y\in[0,1]$, so that
\begin{equation}
    q^y(\rho):=q^y_{\alpha^y,\lambda^y,p^y}(\rho)
\end{equation}
By \eqref{eq:flux_boundary_values}, $q^y$ retains the boundary values $q^y(0)=q^y(\rho_{\max})=0$, and additionally satisfies
\[
q^y(p^y\,\rho_{\max})=0.
\]
The parameter $p^y$ can therefore be interpreted as the density at which the preferred direction of lane changes switches. This choice differs from the reduced two-parameter model proposed in \cite{2D_traffic} for the $y$-direction. We found that the reduced model does not reproduce the fundamental diagrams of the datasets considered here sufficiently well and therefore use the more flexible construction \eqref{eq:flux_param_y}.

Since $q^y$ is defined as a product involving $q_{\alpha,\lambda,p}$, the properties of Proposition~\ref{prop:flux_structure} only partially carry over. The boundary values \eqref{eq:flux_boundary_values} and the $\lambda$-symmetry \eqref{eq:lambda_symmetry} still hold for $q^y$, since both are inherited directly from $q_{\alpha,\lambda,p}$. In contrast, the concavity property \eqref{eq:flux_second_derivative} and the constant-sign property \eqref{eq:constant_sign} do not transfer.

\subsection{Parameter calibration and dimension reduction}
\label{sec:parameter_identification}

To identify the parameter sets $(\alpha^x,\lambda^x,p^x)$ and $(\alpha^y,\lambda^y,p^y)$ of the flux models \eqref{eq:flux_param} and \eqref{eq:flux_param_y}, we formulate independent non-linear least-squares problems for the traffic flux in $x$- and $y$-direction. Therefore, we consider the measured density-flux pairs
\[
(\rho_{k_0},q_{k_0}^z),\qquad k_{0}=0,\ldots,N,\qquad z\in\{x,y\},
\]
where $N=M+1-m$, for each dataset introduced in Section~\ref{sec:data_description}. Once the parameter sets have been determined, the mean speed is recovered from the fundamental relation $q=\rho u$ as
\begin{equation}
\label{eq:velocity_from_flux}
u^z(\rho)=\frac{q^z(\rho)}{\rho},
\qquad \rho>0,
\end{equation}
which extends continuously to $\rho=0$ by $u^z(0):=(q^z)'(0)$, the free-flow speed.

Writing
\begin{equation}
\label{eq:shape_function}
g_{\lambda,p}(\rho) = d_1+(d_2-d_1)\frac{\rho}{\rho_{\max}} - \sqrt{1+\lambda^2\left(\frac{\rho}{\rho_{\max}}-p\right)^2},
\end{equation}
so that
\begin{equation}
\label{eq:shape_function_y}
q^x_{\alpha,\lambda,p}(\rho)=\alpha\,g_{\lambda,p}(\rho), \qquad 
q^y_{\alpha,\lambda,p}(\rho) := \alpha\left(\frac{\rho}{\rho_{\max}}-p\right) g_{\lambda,p}(\rho),
\end{equation}
we define
\begin{equation}
\label{eq:flux_separable}
g^z_{\lambda,p}(\rho):= 
\begin{cases}
    \,g_{\lambda,p}(\rho) & z=x\\
    \left(\frac{\rho}{\rho_{\max}}-p\right) g_{\lambda,p}(\rho) & z=y
\end{cases}.
\end{equation}

The least-squares objective of the parameter calibration problem in direction $z\in\{x,y\}$ therefore takes the form
\begin{equation}
\label{eq:objectives}
J^z(\alpha,\lambda,p) = \frac12\sum_{k_0=0}^{N} \left(\alpha g^z_{\lambda,p}(\rho_{k_0})-q^z_{k_0}\right)^2.
\end{equation}
While \eqref{eq:objectives} is non-linear in $\lambda$ and $p$, it is a strictly convex quadratic function in $\alpha$ whenever the vector
\[
g^z(\lambda,p) = \bigl(g^z_{\lambda,p}(\rho_0),\ldots,g^z_{\lambda,p}(\rho_N)\bigr)^\top
\]
is non-zero. Consequently, one of the three parameters can be eliminated analytically for each $z\in\{x,y\}$. This is a special case of the variable projection principle for separable non-linear least-squares problems \cite{golub_pereyra_1973}. We first state this reduction explicitly for $z=x$. The case $z=y$ then follows with minor adjustments, since $g^y_{\lambda,p}$ vanishes at the additional point $\rho=p\,\rho_{\max}$.

\begin{proposition}[Exact elimination of $\alpha$ for the $x$-direction]
\label{prop:alpha_elimination}
Fix $\lambda\neq0$ and $p\in\mathbb{R}$ and assume that at least one observation satisfies
$0<\rho_{k_0}<\rho_{\max}$.
Then $g_{\lambda,p}(\rho)>0$ for such $\rho_{k_0}$, hence $g^x(\lambda,p)\neq0$, and the unique minimizer of $J^x$ from \eqref{eq:objectives} with respect to $\alpha\in\mathbb{R}$ is
\begin{equation}
    \label{eq:alpha_star}
    \alpha^{\star}(\lambda,p) = \frac{\displaystyle\sum_{k_0=0}^{N}q^x_{k_0}\, g_{\lambda,p}(\rho_{k_0})}{\displaystyle \sum_{k_0=0}^{N}g_{\lambda,p}(\rho_{k_0})^2}.
\end{equation}
The corresponding reduced objective is
\begin{align}
\label{eq:reduced_objective}
    \widehat J^x(\lambda,p) &:=J^x\bigl(\alpha^\star(\lambda,p),\lambda,p\bigr) =
    \frac12\left( \sum_{k_0=0}^{N}(q^x_{k_0})^2 - \frac{\left(\displaystyle\sum_{k_0=0}^{N}q^x_{k_0}\, g_{\lambda,p}(\rho_{k_0})\right)^2} {\displaystyle\sum_{k_0=0}^{N}g_{\lambda,p}(\rho_{k_0})^2} \right).
\end{align}
If $\alpha$ is restricted to a closed interval $I=[\alpha_{\min},\alpha_{\max}]$, the minimizing value is obtained by projection,
\begin{equation}
    \label{eq:alpha_projected}
    \alpha_I^\star(\lambda,p) =
    \begin{cases}
    \alpha_{\min}, & \text{if } \alpha^\star(\lambda,p) < \alpha_{\min},\\[2pt]
    \alpha^\star(\lambda,p), & \text{if } \alpha_{\min} \le \alpha^\star(\lambda,p) \le \alpha_{\max},\\[2pt]
    \alpha_{\max}, & \text{if } \alpha^\star(\lambda,p) > \alpha_{\max}.
    \end{cases}
\end{equation}
\end{proposition}

\begin{proof}
By Proposition~\ref{prop:flux_structure}, $g_{\lambda,p}(\rho)>0$ for every $0<\rho<\rho_{\max}$ whenever $\lambda\neq0$. Hence the assumption on the data implies
\[
A:=\sum_{k_0=0}^{N}g_{\lambda,p}(\rho_{k_0})^2>0.
\]
With $\displaystyle B:=\sum_{k_0=0}^{N}q^x_{k_0}\,g_{\lambda,p}(\rho_{k_0})$ and $\displaystyle C:=\sum_{k_0=0}^{N}(q^x_{k_0})^2$, expanding \eqref{eq:objectives} gives
\[
J^x(\alpha,\lambda,p) = \frac12\left(A\alpha^2-2B\alpha+C\right),
\]
a strictly convex quadratic function in $\alpha$. Its unique global minimum is given by $\alpha^\star=\frac{B}{A}$, which proves \eqref{eq:alpha_star}. Substituting back yields
\[
\widehat J^x(\lambda,p) = \frac12\left(A\left(\frac{B}{A}\right)^2-2B\left(\frac{B}{A}\right)+C\right)
= \frac12\left(C-\frac{B^2}{A}\right),
\]
which gives \eqref{eq:reduced_objective}. Since $J^x(\cdot,\lambda,p)$ is strictly convex with unique minimizer $\alpha^\star$, it is strictly decreasing on $(-\infty,\alpha^\star]$ and strictly increasing on $[\alpha^\star,\infty)$; the minimizer over the closed interval $I$ is therefore $\alpha_{\min}$ if $\alpha^\star<\alpha_{\min}$, $\alpha_{\max}$ if $\alpha^\star>\alpha_{\max}$, and $\alpha^\star$ otherwise.
\end{proof}

\begin{corollary}[Exact elimination of $\alpha$ for the $y$-direction]
\label{cor:alpha_elimination_y}
Fix $\lambda\neq0$ and $p\in[0,1]$ and assume that at least one observation satisfies
\[
0<\rho_{k_0}<\rho_{\max}
\qquad\text{and}\qquad
\rho_{k_0}\neq p\,\rho_{\max}.
\]
Then $g^y(\lambda,p)\neq0$, and the unique minimizer of $J^y$ from \eqref{eq:objectives} with respect to $\alpha\in\mathbb{R}$, the reduced objective $\widehat J^y(\lambda,p)$, and, if $\alpha$ is restricted to $I$, the projected minimizer $\alpha^\star_I(\lambda,p)$, are obtained from \eqref{eq:alpha_star}, \eqref{eq:reduced_objective}, and \eqref{eq:alpha_projected}, respectively, with $g_{\lambda,p}$ and $q^x_{k_0}$ replaced by $\left(\frac{\rho}{\rho_{\max}}-p\right)g_{\lambda,p}$ and $q^y_{k_0}$.
\end{corollary}

\begin{proof}
By Proposition~\ref{prop:flux_structure}, $g_{\lambda,p}(\rho)>0$ for $0<\rho<\rho_{\max}$, so $\left(\frac{\rho}{\rho_{\max}}-p\right)g_{\lambda,p}(\rho)$ vanishes on $(0,\rho_{\max})$ only at $\rho=p\,\rho_{\max}$. The assumption on the data therefore guarantees 
$$A^y:=\sum_{k_0=0}^N \left(\left(\frac{\rho}{\rho_{\max}}-p\right)g_{\lambda,p}(\rho)\right)^2>0.$$
The remainder of the proof of Proposition~\ref{prop:alpha_elimination} depends only on this positivity and applies analogously with $g_{\lambda,p}$ and $q^x_{k_0}$ replaced by $\left(\frac{\rho}{\rho_{\max}}-p\right)g_{\lambda,p}$ and $q^y_{k_0}$.
\end{proof}

Proposition~\ref{prop:alpha_elimination} and Corollary~\ref{cor:alpha_elimination_y} reduce the original three-dimensional calibration problem in each direction $z\in\{x,y\}$ exactly to a two-dimensional problem,
\begin{equation}
\label{eq:reduced_calibration}
\min_{(\lambda,p)\in\mathcal D}\widehat J^z(\lambda,p),
\end{equation}
followed by the reconstruction of $\alpha$ from \eqref{eq:alpha_star} or \eqref{eq:alpha_projected}. Moreover, since $\displaystyle \sum_{k_0=0}^{N}(q^z_{k_0})^2$ does not depend on $\lambda$ and $p$, minimizing \eqref{eq:reduced_objective} is equivalent to maximizing
\begin{equation}
\label{eq:normalized_alignment}
\frac{\left\langle q^z,g^z(\lambda,p)\right\rangle^2}
{\|g^z(\lambda,p)\|_2^2},
\end{equation}
where $q^z=(q^z_0,\ldots,q^z_N)^\top$.

The reduced formulation also yields an elementary existence result. If
\[
\mathcal D=[\lambda_{\min},\lambda_{\max}]\times[0,1],
\qquad 0<\lambda_{\min}<\lambda_{\max},
\]
is compact and the data contain at least one density in $(0,\rho_{\max})$, then, by Proposition~\ref{prop:flux_structure}, $\displaystyle\sum_{k_0=0}^N g_{\lambda,p}(\rho_{k_0})^2$ is continuous and strictly positive on $\mathcal D$, so $\widehat J^x$ is continuous on $\mathcal D$, and, by the Weierstrass extreme value theorem, \eqref{eq:reduced_calibration} admits a global minimizer for $z=x$. For $z=y$, the same argument requires $\sum_{k_0=0}^N \left(\frac{\rho}{\rho_{\max}}-p\right)g_{\lambda,p}(\rho_{k_0})^2$ to be strictly positive over $\mathcal D$. Since $\left(\frac{\rho}{\rho_{\max}}-p\right)g_{\lambda,p}$ vanishes at $\rho=p\,\rho_{\max}$, this holds in particular whenever the data contain at least two distinct densities in $(0,\rho_{\max})$. Under this assumption, $\widehat J^y$ is continuous on $\mathcal D$, and \eqref{eq:reduced_calibration} admits a global minimizer for $z=y$ as well. Notice that the lower bound $\lambda_{\min}>0$ excludes the degenerate case $\lambda=0$, for which $g_{0,p}\equiv0$, $g^y_{0,p}\equiv0$, and $\alpha$ is not identifiable in either direction.

Besides the simplification of the parameter calibration established above, a more subtle form of near-redundancy within the remaining two-dimensional problem \eqref{eq:reduced_calibration} can be analyzed. Substantially different pairs $(\lambda,p)$ can produce nearly identical values of $\widehat J^z$, and hence nearly indistinguishable fitted traffic models.

\begin{proposition}[Parameter degeneracy for the $x$-direction]
\label{prop:large_lambda_asymptotics}
Let $0\leq p\leq1$, $\lambda>0$ and set $r=\frac{\rho}{\rho_{\max}}$. Define the triangular profile
\begin{equation}
    \label{eq:triangular_profile}
    \tau_p(r) = (1-r)p+r(1-p)-|r-p| =
    \begin{cases}
    2r(1-p), & 0\leq r\leq p,\\
    2p(1-r), & p\leq r\leq1.
    \end{cases}
\end{equation}
Then
\begin{equation}
\label{eq:large_lambda_expansion}
    g_{\lambda,p}(\rho)=\lambda\,\tau_p(r)+R_{\lambda,p}(r), \qquad
    |R_{\lambda,p}(r)|\leq1
\end{equation}
uniformly for $r\in[0,1]$. In particular, for $\alpha \neq 0$
\begin{equation}
    \label{eq:large_lambda_uniform_limit}
    \left\| \frac{q^x_{\alpha,\lambda,p}}{\alpha\lambda} - \tau_p\!\left(\frac{\cdot}{\rho_{\max}}\right) \right\|_{L^\infty(0,\rho_{\max})} \leq\frac1\lambda.
\end{equation}
\end{proposition}

\begin{proof}
By computing
\[
\sqrt{1+s^2}-s =\frac{\bigl(\sqrt{1+s^2}-s\bigr)\bigl(\sqrt{1+s^2}+s\bigr)}{\sqrt{1+s^2}+s}
=\frac{1}{\sqrt{1+s^2}+s}
\]
for $s\geq0$, we can write
\[
\sqrt{1+s^2}=s+\varepsilon(s), \qquad \varepsilon(s):=\frac{1}{\sqrt{1+s^2}+s}.
\]
Since $\sqrt{1+s^2}+s\geq\sqrt{1+s^2}\geq1$ for $s\geq0$, we get $0<\varepsilon(s)\leq1$.

Recall from \eqref{eq:shape_function} that
\[
g_{\lambda,p}(\rho) = (1-r)\sqrt{1+(\lambda p)^2} + r\sqrt{1+(\lambda(1-p))^2} - \sqrt{1+\lambda^2(r-p)^2}.
\]
Applying the identity $\sqrt{1+s^2}=s+\varepsilon(s)$ to the three square-root terms, with $s=\lambda p$, $s=\lambda(1-p)$, and $s=\lambda|r-p|$ respectively, gives
\begin{align*}
\sqrt{1+(\lambda p)^2} &= \lambda p+\varepsilon(\lambda p),\\
\sqrt{1+(\lambda(1-p))^2} &= \lambda(1-p)+\varepsilon(\lambda(1-p)),\\
\sqrt{1+\lambda^2(r-p)^2} &= \lambda|r-p|+\varepsilon(\lambda|r-p|).
\end{align*}
Substituting these three expressions into $g_{\lambda,p}(\rho)$ yields
\begin{align*}
g_{\lambda,p}(\rho)
&=(1-r)\bigl[\lambda p+\varepsilon(\lambda p)\bigr]
+r\bigl[\lambda(1-p)+\varepsilon(\lambda(1-p))\bigr]
-\bigl[\lambda|r-p|+\varepsilon(\lambda|r-p|)\bigr]\\
&=\underbrace{\lambda\bigl[(1-r)p+r(1-p)-|r-p|\bigr]}_{=:\,\lambda\tau_p(r)}
+\underbrace{(1-r)\varepsilon(\lambda p)+r\,\varepsilon(\lambda(1-p))-\varepsilon(\lambda|r-p|)}_{=:\,R_{\lambda,p}(r)},
\end{align*}
so we obtain
\[
g_{\lambda,p}(\rho)=\lambda\,\tau_p(r)+R_{\lambda,p}(r).
\]

Since the first two terms of $R_{\lambda,p}$ form a convex combination of the two numbers $\varepsilon(\lambda p),\varepsilon(\lambda(1-p))\in(0,1]$, the convex combination itself lies in $(0,1]$. Hence
\[
-1=0-1< R_{\lambda,p}(r)\leq 1-0=1,
\]
so $|R_{\lambda,p}(r)|\leq1$, uniformly in $r\in[0,1]$ and $\lambda>0$. This proves \eqref{eq:large_lambda_expansion}.

By \eqref{eq:flux_separable}, $q^x_{\alpha,\lambda,p}(\rho)=\alpha\, g_{\lambda,p}(\rho)$, so for $\alpha\neq0$,
\[
\frac{q^x_{\alpha,\lambda,p}(\rho)}{\alpha\lambda} = \frac{g_{\lambda,p}(\rho)}{\lambda} = \tau_p(r)+\frac{R_{\lambda,p}(r)}{\lambda}
\]
by \eqref{eq:large_lambda_expansion}. Taking absolute values and using $|R_{\lambda,p}(r)|\le1$ gives
\[
\left|\frac{q^x_{\alpha,\lambda,p}(\rho)}{\alpha\lambda}-\tau_p(r)\right|\leq\frac1\lambda
\]
for every $\rho\in(0,\rho_{\max})$. Taking the supremum over $\rho$ yields \eqref{eq:large_lambda_uniform_limit}.
\end{proof}

\begin{corollary}[Parameter degeneracy for the $y$-direction]
\label{cor:large_lambda_asymptotics_y}
Let $0\leq p\leq1$, $\lambda>0$ and set $r=\frac{\rho}{\rho_{\max}}$. Define the signed triangular profile
\begin{equation}
\label{eq:triangular_profile_y}
\tau_p^y(r):=(r-p)\,\tau_p(r).
\end{equation}
Then
\begin{equation}
\label{eq:large_lambda_expansion_y}
\left(\frac{\rho}{\rho_{\max}}-p\right)g_{\lambda,p}(\rho) = \lambda\,\tau_p^y(r) + R^y_{\lambda,p}(r),
\qquad
|R^y_{\lambda,p}(r)|\leq1
\end{equation}
uniformly for $r\in[0,1]$. In particular, for $\alpha\neq0$,
\begin{equation}
\label{eq:large_lambda_uniform_limit_y}
\left\|
\frac{q^y_{\alpha,\lambda,p}}{\alpha\lambda}
-\tau_p^y\!\left(\frac{\cdot}{\rho_{\max}}\right)
\right\|_{L^\infty(0,\rho_{\max})}
\leq\frac1\lambda.
\end{equation}
\end{corollary}

\begin{proof}
Multiplying \eqref{eq:large_lambda_expansion} by $(r-p)$ gives
\[
\left(\frac{\rho}{\rho_{\max}}-p\right)g_{\lambda,p}(\rho)=(r-p)\,g_{\lambda,p}(\rho)=\lambda(r-p)\tau_p(r)+(r-p)R_{\lambda,p}(r),
\]
which is \eqref{eq:large_lambda_expansion_y} with $R^y_{\lambda,p}(r):=(r-p)R_{\lambda,p}(r)$. Since $r,p\in[0,1]$, we have $|r-p|\leq1$, so $|R^y_{\lambda,p}(r)|\leq|R_{\lambda,p}(r)|\leq1$ by \eqref{eq:large_lambda_expansion}. Estimate \eqref{eq:large_lambda_uniform_limit_y} then follows exactly as in the proof of Proposition~\ref{prop:large_lambda_asymptotics}, using \eqref{eq:flux_separable} for $z=y$.
\end{proof}

Proposition~\ref{prop:large_lambda_asymptotics} and Corollary~\ref{cor:large_lambda_asymptotics_y} yield for large $\lambda$ that
\begin{equation}
\label{eq:large_lambda_flux_interpretation}
q^x_{\alpha,\lambda,p}(\rho)
\approx
(\alpha\lambda)\,
\tau_p\!\left(\frac{\rho}{\rho_{\max}}\right), \qquad q^y_{\alpha,\lambda,p}(\rho)\approx(\alpha\lambda)\,\tau_p^y\!\left(\frac{\rho}{\rho_{\max}}\right).
\end{equation}
Thus $p$ determines the location of the turning point of the curve, whereas the slope depends mostly on the product $\alpha\lambda$. Separate changes in $\alpha$ and $\lambda$ that approximately preserve this product therefore lead to nearly identical flux functions. In addition, it confirms the interpretation of $p^y$ given in Section~\ref{sec:data_fitting}. Since $\tau_p^y$ changes sign exactly at $r=p$, the limiting $y$-flux switches sign at $\rho=p\,\rho_{\max}$.

\begin{remark}
Complementary to Proposition~\ref{prop:large_lambda_asymptotics} and Corollary~\ref{cor:large_lambda_asymptotics_y}, a Taylor expansion of $\sqrt{1+s^2}=1+\frac{s^2}{2}+O(s^4)$ in the three square-root terms of $g_{\lambda,p}$, analogous to the proof of Proposition~\ref{prop:large_lambda_asymptotics}, yields for $\lambda\to0$
\[
g_{\lambda,p}(\rho)=\frac{\lambda^2}{2}\,r(1-r)+O(\lambda^4),
\qquad r=\frac{\rho}{\rho_{\max}},
\]
uniformly in $r,p\in[0,1]$. Hence
\[
q^x_{\alpha,\lambda,p}(\rho)\approx\frac{\alpha\lambda^2}{2}\,r(1-r),
\qquad
q^y_{\alpha,\lambda,p}(\rho)\approx\frac{\alpha\lambda^2}{2}\,(r-p)\,r(1-r).
\]
In the $x$-direction, the flux depends on $\alpha$ and $\lambda$ only through $\alpha\lambda^2$ and no longer on $p$. In the $y$-direction, $p$ remains. Together with Proposition~\ref{prop:large_lambda_asymptotics}, this shows that the parameters are only weakly identifiable at both ends of the $\lambda$-range.
\end{remark}

\medskip

The simplified optimization problem \eqref{eq:reduced_calibration} can now be solved independently for $z\in\{x,y\}$ and each dataset $(\rho_{k_0},q_{k_0}^z)$. Since the objective functions $\widehat J^z$ are non-linear and non-convex, efficient numerical optimization techniques that are capable of escaping local minima are required.

In the following section, we introduce and compare several optimization algorithms for solving these parameter calibration problems. The methods range from classical gradient-based approaches to modern derivative-free optimization algorithms.

\section{Optimization algorithms}
\label{sec:optimization}

The preceding analysis leads to two equivalent formulations of the calibration problem. The \emph{full formulation} minimizes the three-parameter least-squares functional $J^z(\alpha,\lambda,p)$ from \eqref{eq:objectives}, whereas the \emph{reduced formulation} minimizes $\widehat J^z(\lambda,p)$ from \eqref{eq:reduced_objective} and reconstructs the optimal $\alpha$ by \eqref{eq:alpha_star}.

The optimization algorithms considered below can be applied to either formulation. We compare two classical local optimization methods, namely \textit{Nelder--Mead} and \textit{L-BFGS-B}, one evolutionary population-based algorithm, \textit{Differential Evolution}, and two interacting-particle methods, \textit{Consensus-Based Optimization} (CBO) and its deterministic variant \textit{QMC$^2$BO}. We present the algorithms in a dimension-independent form so that the same implementation can be used with $D=3$ for the full problem and $D=2$ for the reduced problem. For each method, we briefly summarize the underlying algorithmic idea before describing the Python implementation and the parameter settings used throughout the numerical experiments. The main characteristics of all algorithms are finally summarized in Table \ref{tab:optimization_methods}.

\subsection{General optimization problem and calibration formulations}
\label{sec:general_optimization_problem}

Although all five algorithms are based on different ideas, they all aim to solve optimization problems of the form
\begin{equation}
\min_{\theta\in\mathcal{K}} J_{\mathrm{opt}}(\theta),
\label{eq:general_optimization_problem}
\end{equation}
where $\theta\in\mathbb{R}^{D}$ denotes the optimization variables, $J_{\mathrm{opt}}:\mathcal{K}\to\mathbb{R}$ the objective function, and $\mathcal{K}\subset\mathbb{R}^{D}$ a compact feasible set. In our setting, there are two choices. For readability, we write $J$ instead of $J_{\mathrm{opt}}$ in the algorithmic formulas below.

\begin{description}

\item[\textbf{Full formulation.}] The original formulation uses
\begin{equation}
\label{eq:full_optimization_variables}
\theta=(\alpha,\lambda,p)\in\mathbb{R}^3,
\qquad
J_{\mathrm{opt}}=J^z,
\end{equation}
with $J^z$ defined in \eqref{eq:objectives}. In view of Proposition~\ref{prop:flux_structure}, i.e. the symmetry \eqref{eq:lambda_symmetry}, we restrict the feasible set to
\begin{equation}
\label{eq:full_admissible_box}
\mathcal K = [-500,500]\times[0,500] \times [0,1].
\end{equation}


\item[\textbf{Reduced formulation.}]
Using Proposition~\ref{prop:alpha_elimination} and its projected version \eqref{eq:alpha_projected} with $I=[-500,500]$, the same least-squares problem can be written with only two non-linear variables,
\begin{equation}
\label{eq:reduced_optimization_variables}
\theta=(\lambda,p)\in\mathbb{R}^2, \qquad J_{\mathrm{opt}}=\widehat J^z, \qquad \mathcal K=[0,500] \times [0,1].
\end{equation}
Once a minimizer $(\lambda^\star,p^\star)$ has been found, we recover $\alpha$ exactly as  $\alpha^\star=\alpha_I^\star(\lambda^\star,p^\star)$ from \eqref{eq:alpha_projected}. 

\end{description}

For all algorithms, we denote the approximation of the optimizer after the $k$-th iteration by $\theta_k$. The corresponding update rule can be written abstractly as
\begin{equation}
\theta_{k+1}=\Phi_k(\theta_k),
\label{eq:general_iteration}
\end{equation}
where the mapping $\Phi_k$ depends on the specific optimization algorithm. Some algorithms evolve a set of candidate solutions rather than a single point. In these cases, we write
\[
\Theta_k=\left(\theta_k^{1},\ldots,\theta_k^{N}\right)\in\mathbb{R}^{D\times N},
\]
where $N$ is the number of possible solutions and $\theta_k^{n}$ denotes the $n$-th one. For particle-based methods, these correspond to the particle positions, whereas for the Nelder-Mead algorithm, they represent the vertices of the current simplex. The approximation of the optimizer, which is denoted by $\theta_k$, is obtained from the current set of possible solutions. For particle-based methods, it is the consensus point, for the Nelder-Mead algorithm it is the best vertex.

We note that we will restrict each optimization method to $10^5$ function evaluations for the full calibration problem \eqref{eq:full_optimization_variables} and $10^4$ function evaluations for the reduced one \eqref{eq:reduced_optimization_variables}. This secures the comparability of each method with respect to computational effort.

\subsection{Nelder-Mead}
\label{sec:nelder_mead}

By successively improving a set of points through geometric transformations, the Nelder-Mead algorithm \cite{ nelder_mead_python_source,nelder_mead_origin} is a derivative-free method which relies solely on evaluations of the objective function. The algorithm maintains a simplex consisting of $D+1$ vertices $\theta_k^1,\ldots,\theta_k^{D+1}\in\mathbb{R}^{D}$. At the beginning of every iteration, the vertices are sorted according to their objective values, i.e.
$$
J(\theta_k^1)\leq J(\theta_k^2)\leq\cdots\leq J(\theta_k^{D+1}),
$$
such that $\theta_k^1$ denotes the current best point and $\theta_k^{D+1}$ the current worst point. The centroid of all vertices except the worst one is then computed as
$$
\bar{\theta}_k = \frac1D\sum_{d=1}^{D}\theta_k^{d}.
$$
Based on this centroid, the algorithm replaces the worst vertex with another point. Depending on the objective value, one of the following points is accepted.

\begin{description}

\item[\textbf{Reflection.}]
The first candidate point is obtained by reflecting the worst vertex at the centroid,
$$
\theta_r=\bar{\theta}_k+\alpha(\bar{\theta}_k-\theta_k^{D+1}),
$$
where $\alpha>0$ denotes the reflection parameter. If now $J(\theta_k^{1})\leq J(\theta_r)<J(\theta_k^{D})$, the reflected point replaces the worst vertex, i.e $\theta_{k+1}^{D+1} = \theta_r$.

\item[\textbf{Expansion.}]
If $J(\theta_r)<J(\theta_k^{1})$, the algorithm continues in the same direction by computing
$$
\theta_e=\bar{\theta}_k+\beta(\theta_r-\bar{\theta}_k),
$$
where $\beta>1$ is the expansion parameter. If $J(\theta_e)<J(\theta_r)$, the expanded point replaces the worst vertex, i.e. $\theta_{k+1}^{D+1} = \theta_e$. Otherwise, the reflected point replaces it, i.e. $\theta_{k+1}^{D+1} = \theta_r$.

\item[\textbf{Outside contraction.}]
If $J(\theta_k^{D})\leq J(\theta_r)<J(\theta_k^{D+1})$, an outside contraction is performed,
$$
\theta_{oc} = \bar{\theta}_k + \gamma(\theta_r-\bar{\theta}_k),
$$
where $0<\gamma<1$. This reduces the step length but keeps the same direction. If $J(\theta_{oc})\leq J(\theta_r)$, the outside contraction point replaces the current worst vertex, ie. $\theta_{k+1}^{D+1} = \theta_{oc}$. Otherwise, the algorithm proceeds with the shrink step.

\item[\textbf{Inside contraction.}]
If $J(\theta_r)\geq J(\theta_k^{D+1})$, the simplex is contracted towards its centroid by computing
$$
\theta_{ic} = \bar{\theta}_k - \gamma(\bar{\theta}_k-\theta_k^{(D+1)}).
$$
and replacing the current worst vertex with $\theta_{ic}$ if $J(\theta_{ic})<J(\theta_k^{D+1})$, i.e. $\theta_{k+1}^{D+1} = \theta_{ic}$. Otherwise, the algorithm proceeds with the shrink step.

\item[\textbf{Shrink.}]
The complete simplex is shrunk towards the current best vertex,
$$
\theta_{k+1}^{d} = \theta_k^{1} + \delta(\theta_k^{d}-\theta_k^{1}), \qquad d=2,\ldots,D+1,
$$
where $0<\delta<1$. This reduces the search region and allows the algorithm to continue from a smaller one.

\end{description}

\begin{figure}[ht]
\centering
\begin{tikzpicture}[
    scale=0.75,
    every node/.style={font=\small},
    vertex/.style={circle,fill=black,inner sep=1.8pt},
    newvertex/.style={circle,fill=red,inner sep=2.0pt},
    centroid/.style={circle,fill=blue,inner sep=1.8pt},
    simplex/.style={thick,gray},
    helper/.style={dashed,gray},
    arrow/.style={->,thick,gray}
]

\begin{scope}[shift={(0,0)}]
\node at (1.5,3.0) {\textbf{Reflection}};
\coordinate (A) at (0,0);
\coordinate (B) at (2,0.2);
\coordinate (C) at (1.5,2);
\draw[simplex] (A)--(B)--(C)--cycle;
\node[vertex] at (A){};
\node[vertex] at (B){};
\node[vertex] at (C){};
\coordinate (M) at (1,0.1);
\node[centroid] at (M){};
\coordinate (R) at (0.5,-1.8);
\node[newvertex] at (R){};
\draw[->,gray,thick,dashed] (C)--(R);
\end{scope}

\begin{scope}[shift={(4.2,0)}]
\node at (1.5,3.0) {\textbf{Expansion}};
\coordinate (A) at (0,0);
\coordinate (B) at (2,0.2);
\coordinate (C) at (1.5,2);
\draw[simplex] (A)--(B)--(C)--cycle;
\node[vertex] at (A){};
\node[vertex] at (B){};
\node[vertex] at (C){};
\coordinate (M) at (1,0.1);
\node[centroid] at (M){};
\coordinate (R) at (0.5,-1.8);
\coordinate (E) at (0.2,-3);
\node[newvertex] at (E){};
\draw[->,gray,thick,dashed] (C)--(E);
\end{scope}

\begin{scope}[shift={(8.4,0)}]
\node at (1.5,3.0) {\textbf{Outside Contr.}};
\coordinate (A) at (0,0);
\coordinate (B) at (2,0.2);
\coordinate (C) at (1.5,2);
\draw[simplex] (A)--(B)--(C)--cycle;
\node[vertex] at (A){};
\node[vertex] at (B){};
\node[vertex] at (C){};
\coordinate (M) at (1,0.1);
\node[centroid] at (M){};
\coordinate (OC) at (0.76,-0.8);
\node[newvertex] at (OC){};
\draw[->,gray,thick,dashed] (C)--(OC);
\end{scope}

\begin{scope}[shift={(12.6,0)}]
\node at (1.5,3.0) {\textbf{Inside Contr.}};
\coordinate (A) at (0,0);
\coordinate (B) at (2,0.2);
\coordinate (C) at (1.5,2);
\draw[simplex] (A)--(B)--(C)--cycle;
\node[vertex] at (A){};
\node[vertex] at (B){};
\node[vertex] at (C){};
\coordinate (M) at (1,0.1);
\node[centroid] at (M){};
\coordinate (IC) at (1.25,1.1);
\node[newvertex] at (IC){};
\draw[->,gray,thick,dashed] (C)--(IC);
\end{scope}

\begin{scope}[shift={(16.8,0)}]
\node at (1.5,3.0) {\textbf{Shrink}};
\coordinate (A) at (0,0);
\coordinate (B) at (2,0.2);
\coordinate (C) at (1.5,2);
\draw[simplex] (A)--(B)--(C)--cycle;
\node[vertex] at (A){};
\node[vertex] at (B){};
\node[vertex] at (C){};
\coordinate (A2) at (0,0);
\coordinate (B2) at (0.9,0.1);
\coordinate (C2) at (0.7,0.9);
\node[newvertex] at (B2){};
\node[newvertex] at (C2){};
\draw[gray,dashed,->] (B) -- (B2);
\draw[gray,dashed,->] (C) -- (C2);
\end{scope}

\end{tikzpicture}

\caption{Geometric interpretation of all possible Nelder-Mead replacements in two dimensions. The black points denote the simplex vertices at iteration $k$, with the upper vertex representing the current worst point and the lower-left vertex the current best point. The blue point marks the centroid, while the dashed lines indicate the possible movements of the new vertex generated by reflection, expansion, or contraction. The red point(s) represent the vertex (vertices) added to the simplex in iteration $k+1$.}
\label{fig:nelder_mead_geometry}
\end{figure}
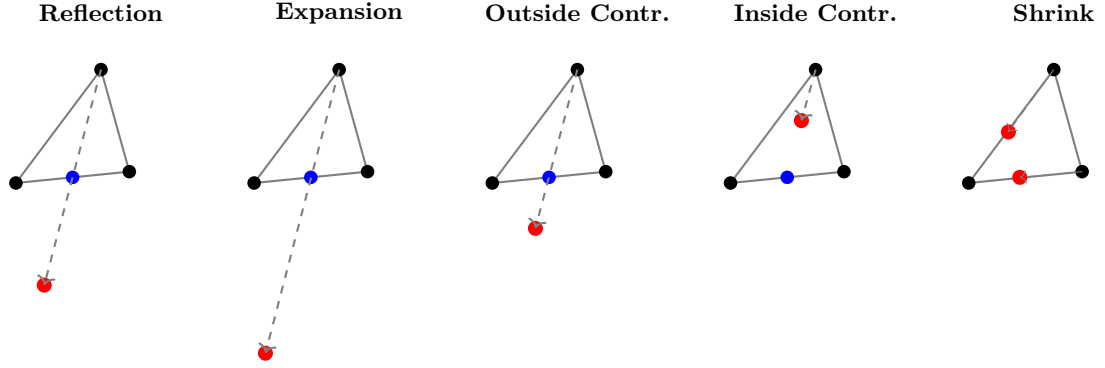

Figure~\ref{fig:nelder_mead_geometry} illustrates the geometric interpretation of each possible replacement for the worst vertex in the two-dimensional case.

Although the Nelder-Mead method lacks a complete convergence theory and may become inefficient for high-dimensional optimization problems \cite{nelder_mead_low_dim_2,nelder_mead_high_dim_2,nelder_mead_high_dim_1}, it has proven to be effective for low-dimensional problems \cite{nelder_mead_low_dim_1, nelder_mead_low_dim_2}. Since the optimization problems considered in this work involve only three variables, the method remains promising.

We employ the implementation provided by the Python package \texttt{scipy.optimize.minimize} using the optimization method \texttt{"Nelder-Mead"} \cite{nelder_mead_python_source,scipy2020}. We adopt the standard parameter choices $(\alpha,\beta,\gamma,\delta) = \left(1,2,\frac{1}{2},\frac{1}{2}\right)$ for the Nelder-Mead method. The choice of the initial simplex is discussed separately in Section~\ref{sec:numerical_experiments}. The algorithm is stopped when the objective function of all vertices are close. More precisely, the iteration is terminated if
$$
\max_{d=2,\ldots D+1} |J(\theta_k^d)-J(\theta_k^{1})| \leq  10^{-6}.
$$

\subsection{L-BFGS-B}
\label{sec:lbfgsb}

L-BFGS-B belongs to the class of gradient-based optimization methods. It is a quasi-Newton method which approximates second-order information instead of computing the Hessian matrix. The method is based on the limited-memory BFGS (L-BFGS) algorithm \cite[Chapter 7.2]{nocedal2006numerical}, extended by box constraints, which ensure that every iterate remains within the feasible region $\mathcal{K}$.

For unconstrained optimization problems, Newton's second order method updates the current 
approximation of the minimizer according to
\begin{equation}
\theta_{k+1} = \theta_k - \left(\nabla^2 J(\theta_k)\right)^{-1} \nabla J(\theta_k),
\end{equation}
where $\nabla^2 J(\theta_k)$ denotes the Hessian matrix of the objective function. Under standard regularity assumptions and provided that the initial iterate is sufficiently close to a global minimizer with a positive definite Hessian, Newton's method converges quadratically \cite[Chapter 3.3]{nocedal2006numerical}. However, the explicit computation and storage of the Hessian matrix becomes inefficient for large dimensions, since the Hessian contains $D^2$ entries and solving the Newton system may require $\mathcal{O}(D^3)$ operations.

To overcome this, the BFGS method replaces the inverse Hessian by a sequence of symmetric positive definite matrices $H_k$ satisfying
\begin{equation*}
H_k \approx \left(\nabla^2 J(\theta_k)\right)^{-1}.
\end{equation*}
The iteration step is then given by
\begin{equation}
\theta_{k+1} = \theta_k - \alpha_k H_k \nabla J(\theta_k),
\end{equation}
where $\alpha_k$ denotes the step-size. The approximation $H_k$ is updated using the differences between two consecutive iterates and their gradients,
\begin{equation*}
s_k=\theta_{k+1}-\theta_k, \qquad
y_k=\nabla J(\theta_{k+1})-\nabla J(\theta_k),
\end{equation*}
and setting
\begin{equation}
H_{k+1} = \left(I - \frac{s_{k}y_{k}^T}{y_{k}^T s_{k}}\right) H_{k} \left(I - \frac{y_{k}s_{k}^T}{y_{k}^T s_{k}}\right) + \frac{s_{k}s_{k}^T}{y_{k}^T s_{k}}
\end{equation}

Although BFGS avoids the computation of the Hessian, storing the approximation $H_k$ still requires $\mathcal{O}(D^2)$ memory. The limited-memory BFGS method avoids this storage by not explicitly computing $H_k$. Instead, only the most recent $m$ pairs $(s_i,y_i)$ for $i=k-m,\ldots,k-1$ are stored. These pairs contain the most relevant curvature information and allow the approximation $H_k$ to be reconstructed. After each iteration, 
the oldest pair is discarded and replaced by the newest pair $(s_{k+1},y_{k+1})$. Typical choices for the memory parameter are $m\in\{3,\ldots,20\}$.

We employ the implementation provided by the Python package \texttt{scipy.optimize.minimize} using the method \texttt{"L-BFGS-B"} \cite{scipy2020}. The choice of the initial iterate $\theta_0$ is discussed separately in Section~\ref{sec:numerical_experiments}. The step-size $\alpha_k$ is determined by a line search satisfying the strong Wolfe conditions \cite[Chapter~3.1]{nocedal2006numerical}. The memory parameter $m=10$ is used. The algorithm terminates once the relative change of the objective values satisfies
$$
\frac{ J(\theta_k)-J(\theta_{k+1}) }{ \max\!\left( |J(\theta_k)|,\, |J(\theta_{k+1})|,\, 1 \right) } \leq 10^{-6}.
$$

\subsection{Differential evolution}
\label{sec:differential_evolution}

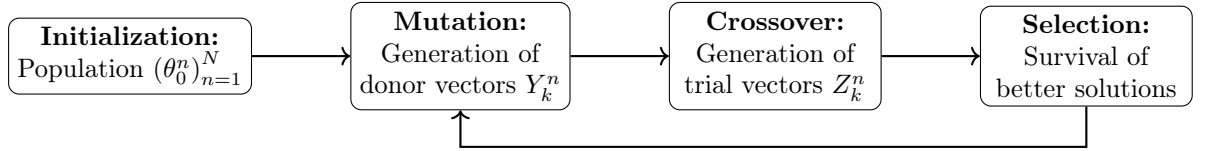
\begin{figure}[ht]
\centering
\begin{tikzpicture}[
    node distance=1.3cm,
    block/.style={draw, rounded corners, minimum width=2.8cm,
    minimum height=0.8cm, align=center},
    arrow/.style={->, thick}
]

\node[block] (init) {\textbf{Initialization:}\\
Population $\left(\theta_0^n\right)_{n=1}^{N}$};
\node[block, right=of init] (mutation) {\textbf{Mutation:} \\ Generation
of \\ donor vectors $Y_k^n$};
\node[block, right=of mutation] (cross) {\textbf{Crossover:} \\ Generation of \\ trial vectors $Z_k^n$};
\node[block, right=of cross] (selection) {\textbf{Selection:} \\ Survival
of \\ better solutions};

\draw[arrow] (init) -- (mutation);
\draw[arrow] (mutation) -- (cross);
\draw[arrow] (cross) -- (selection);
\draw[arrow] (selection) |- ++(0,-1.2) -| 
(node cs:name=mutation,anchor=south);

\end{tikzpicture}
\caption{Framework of the four phases of the Differential Evolution method.}
\label{fig:de_framework}
\end{figure}

Differential Evolution (DE) \cite{differential_evolution_origin} is a population-based optimization method that evolves a set of possible solutions simultaneously and does not require gradient information. More precisely, it is inspired by Darwin's theory of natural evolution and mimics mutation, crossover, and selection to improve a set of possible solutions over several generations. Through exchange between these solutions, the population aims to gradually evolve towards lower objective values. The framework of Differential Evolution consists of four main phases: initialization, mutation, crossover, and selection \cite{differential_evolution_overview}. An overview of the iterative process is illustrated in Figure~\ref{fig:de_framework}.

\begin{description}

\item[\textbf{Initialization.}]
During the initialization phase, a population of $N$ possible solutions is generated. Each solution is represented by a $D$-dimensional vector $\theta_k^{n}$ for $n=1,\ldots,N$, where $k$ denotes the current generation. At $k=0$, each component is uniformly sampled within the feasible parameter set $\mathcal{K}$.

\item[\textbf{Mutation.}]
During the mutation phase, new candidate solutions, called donor vectors, are generated by combining existing solutions. Different mutation strategies exist depending on the number of differences and the selection of the base solution \cite{differential_evolution_overview}. In this work, we employ the \texttt{DE/best/1} strategy, where the donor vector is based on the current best solution and mutated by a single difference
\begin{equation}
Y_k^n = \theta_k^{\mathrm{best}} + \alpha \left( \theta_k^{r_1(n)}-\theta_k^{r_2(n)} \right), \qquad n=1,\ldots N,
\end{equation}
where $r_1(n)\neq r_2(n)\neq n$ denote randomly selected indices and $\alpha>0$ is the mutation factor.

\item[\textbf{Crossover.}]
During the crossover phase the donor vectors $Y_k^n$ are crossed with the original vectors $\theta_k^n$ in a probabilistic manner to generate trial vectors $Z_k^n$. In this work, we use the binomial crossover strategy, where, for each component, the donor vectors are accepted with probability given by the crossover rate $cr\in[0,1]$,
\begin{equation}
Z_k^{n,d} = 
\begin{cases} 
Y_k^{n,d}, & r_{n,d}\leq cr \text{ or }d=d_{\mathrm{rand}}, \\ 
\theta_k^{n,d}, & \text{otherwise}, 
\end{cases} \qquad n=1,\ldots,N, \quad d=1,\ldots,D
\end{equation}
where $r_{n,d}\sim \mathcal{U}(0,1)$ and $d_{\mathrm{rand}}\in \{1,\ldots,D\}$ is a randomly selected index. Another possible crossover strategy is the exponential crossover. For further details we refer to \cite{differential_evolution_overview}.

\item[\textbf{Selection.}]
Finally, during the selection phase it is determined which solution survive to the next 
generation. The trial vector $Z_k^{n}$ replaces the original vector $\theta_k^n$ only if it yields a lower objective value,
\begin{equation}
\theta_{k+1}^{n} =
\begin{cases}
Z_k^n, & J(Z_k^n)\leq J(\theta_k^n), \\
\theta_k^n, & \text{otherwise}.
\end{cases}
\end{equation}
This selection ensures that the population cannot get worse between generations.

\end{description}

After the one-time initialization, the mutation, crossover, and selection phases are repeated until the termination criterion is met.

We employ the Python package \texttt{scipy.optimize.differential\_evolution} 
\cite{scipy2020}. The crossover rate is set to $cr=0.7$ and the mutation factor $\alpha$ is selected uniformly, independently in every generation from the interval $[0.5,1)$. The population size is determined by $N=15 \times D$. The algorithm is stopped when the objective function of all possible solutions are close. More precisely, the iteration is terminated if
$$
\frac{\operatorname{std}\!\left( J(\theta_k^1),\ldots,J(\theta_k^{N}) \right)}{\left| \frac1N \sum_{n=1}^{N} J(\theta_k^n) \right|}
\leq 10^{-6}.
$$

\subsection{Consensus-based optimization}
\label{sec:cbo}

Consensus-Based Optimization (CBO) extends the population-based optimization by modeling the possible solutions as interacting particles evolving according to a stochastic differential equation. It is a gradient-free method for finding global minimizers of possibly non-convex objective functions \cite{CBO_origin,carrillo_analytical_framework_CBO}.

The algorithm governs a finite set of $N$ particles $\theta_t^1,\ldots,\theta_t^N$, which represent possible solutions of the optimization problem. The particles are pushed towards a weighted consensus point
\begin{equation}
\nu_J^\alpha(\Theta_t) = \sum_{m=1}^{N} \theta_t^m \frac{\exp(-\alpha J(\theta_t^m))} {\sum_{k=1}^{N}\exp(-\alpha J(\theta_t^k))},
\label{eq:cbo_consensus}
\end{equation}
where $\alpha>0$ controls the concentration around particles with small objective values.

The dynamics of each particle are governed by the stochastic differential equation
\begin{equation}
d\theta_t^n = -\lambda \left( \theta_t^n-\nu_J^\alpha(\Theta_t) \right)dt + \sigma \left| \theta_t^n - \nu_J^\alpha(\Theta_t) \right| dW_t^n , \qquad t\geq0
\label{eq:cbo_system}
\end{equation}
for $n=1,\ldots,N$. The first term, controlled by $\lambda>0$, is referred to as the \textit{exploitation term} and drives the particles towards consensus. The second term, controlled by $\sigma>0$, is the \textit{exploration term}, which introduces stochastic exploration through independent $D$-dimensional Brownian motions $W_t^n$. 

The continuous-time dynamics \eqref{eq:cbo_system} are usually discretized by the Euler-Maruyama method \citenote{\text{p. 340}}{Num_SDE}. This yields the discrete particle system $\Theta_0,\Theta_1,\Theta_2,\ldots$, introduced in Section~\ref{sec:general_optimization_problem}, where $\Theta_k=(\theta_k^1,\ldots,\theta_k^N)$ denotes the particle positions after the $k$-th iteration. 

The convergence properties of CBO and its mean-field limit have been extensively 
studied, including convergence towards global minimizers and stability properties 
of the corresponding stochastic particle system and its discretizations 
\cite{bonandin2025,carrillo_analytical_framework_CBO,
anisotropic_convergence_CBO,fornasier_CBO_converge_globally,finit_N_CBO_Ha,CBO_origin,riedl_CBO_behaves_SGD_like}.

We employ the Python CBO framework introduced in \cite{CBO_in_python}. Unless stated otherwise, the numerical experiments use a population of $N=2^{10}$ particles for the full calibration problem \eqref{eq:full_optimization_variables} and $N=2^{7}$ particles for the reduced one \eqref{eq:reduced_optimization_variables}, since fewer particles are necessary to solve the smaller two-dimensional optimization problem. The initial particles are sampled independently and uniformly from the feasible parameter box $\mathcal K$, so that the same implementation applies to the full and reduced calibration problems. The consensus parameters are fixed to $\lambda=1$, $\sigma=1$ and $\alpha=100$; these are algorithmic CBO parameters and should not be confused with the parameters of the traffic flux. We use the Euler--Maruyama discretization with time step-size $\Delta t=0.05$ to discretize the system.

In contrast to the other optimization algorithms, CBO is typically stopped after reaching a fixed time horizon. Throughout this work, we perform $96$ Euler-Maruyama steps for the full calibration problem \eqref{eq:full_optimization_variables} and $77$ for the reduced one \eqref{eq:reduced_optimization_variables}. We choose this number of steps to maximize the number of function evaluations $\texttt{feval}$, which for CBO is given by
\begin{equation*}
\texttt{feval} = N \times (\#\text{Euler-Maruyama steps}+1),
\end{equation*}
subject to the constraint that $\texttt{feval}$ does not exceed $10^5$ (respectively $10^4$ for the reduced system).

To ensure a fair comparison with the remaining optimization methods, the final optimizer approximation is projected onto $\mathcal{K}$. Since the stochastic exploration term may drive particles outside this domain, this might be necessary.

\subsection{QMC$^2$BO}
\label{sec:qmc2bo}

QMC$^2$BO is a deterministic variant of CBO which is obtained by removing the exploration term \cite{QMC2BO}. This leads to the following interacting particle system:
\begin{equation}
\frac{d \theta_t^n}{dt} = -\lambda \left( \theta_t^n-\nu_J^\alpha(\Theta_t) \right), \qquad n=1,\ldots,N , \quad t\geq 0
\label{eq:qmc2bo_system}
\end{equation}
where the consensus point $\nu_J^\alpha$ is defined as in \eqref{eq:cbo_consensus}.

Each particle is attracted by the weighted average of all particles, where particles with smaller objective values contribute more to the average. Consequently, agents still interact through exponentially weighted averages biased toward lower function values.

Additionally, the QMC$^2$BO method uses a deterministic choice of initial particle positions. More precisely, let the feasible set be the box
\[
\mathcal K=\prod_{d=1}^{D}[a_d,b_d].
\]
The particles are initialized as
\begin{equation*}
\theta_0^n = \mathfrak{s}(\boldsymbol{p}^n), \qquad n=1,\ldots,N,
\end{equation*}
where $\boldsymbol{p}^1,\ldots,\boldsymbol{p}^N\in[0,1)^D$ are low-discrepancy points and the componentwise affine mapping
\begin{equation}\label{eq:linear_transformation_construction}
\mathfrak{s}_d(u)=a_d+(b_d-a_d)u_d,
\qquad d=1,\ldots,D,
\end{equation}
rescales the unit cube to $\mathcal K$. In this way the deterministic initialization can be used without modification for both the three-dimensional full calibration and the two-dimensional reduced calibration.

We employ the Python CBO framework introduced in \cite{CBO_in_python} with $\sigma=0$ and the corresponding deterministic initialization. Unless stated otherwise, the point set consists of Sobol points \cite{Sobol_points}, generated with the \texttt{scipy.stats.qmc} module \cite{scipy2020} and mapped to the respective feasible parameter box by \eqref{eq:linear_transformation_construction}.

\subsection{Comparison of the optimization algorithms}

Before turning to the numerical experiments, we briefly summarize the main characteristics of the five considered optimization algorithms. Although all methods solve the same parameter identification problem introduced in Section~\ref{sec:parameter_identification}, they differ in their origin, in their optimization strategy and in their robustness with respect to objective functions. Table~\ref{tab:optimization_methods} compares some of their properties.

\begin{table}[ht]
\caption{Comparison of the five optimization algorithms.}
\centering
\renewcommand{\arraystretch}{1.25}

\begin{tabular}{p{3.2cm}ccccc}
\toprule
\textbf{Property}
&
\textbf{Nelder-Mead}
&
\textbf{L-BFGS-B}
&
\textbf{DE}
&
\textbf{CBO}
&
\textbf{QMC$^2$BO}
\\
\midrule

Derivative-free
&
\checkmark
&
\texttimes
&
\checkmark
&
\checkmark
&
\checkmark
\\

Particle-based
&
\texttimes
&
\texttimes
&
\checkmark
&
\checkmark
&
\checkmark
\\

Global optimization
&
\texttimes
&
\texttimes
&
\checkmark
&
\checkmark
&
\checkmark
\\

Deterministic
&
\checkmark
&
\checkmark
&
\texttimes
&
\texttimes
&
\checkmark
\\

Initialization
&
predefined
&
predefined
&
uniformly
&
uniformly
&
Sobol points
\\

\bottomrule
\end{tabular}

\label{tab:optimization_methods}
\end{table}

\section{Numerical experiments}
\label{sec:numerical_experiments}

We apply the five optimization algorithms from Section~\ref{sec:optimization} to the calibration problems associated with the three traffic datasets from Section~\ref{sec:data}. We first compare and discuss the results for the full calibration \eqref{eq:full_optimization_variables} and the reduced calibration \eqref{eq:reduced_optimization_variables} in case of the Off-ramp dataset. We then compare the full calibration of the Off-ramp dataset to the one of the highD dataset and the Heidelberg sensor dataset. Finally, we investigate sensitivity with respect to $\lambda$ and $p$ and robustness with respect to stochasticity and initialization. 

For each dataset, the fundamental diagrams and speed-density diagrams from Figure \ref{fig:fundamental_diagrams_offramp}, \ref{fig:fundamental_diagrams_highD} and \ref{fig:fundamental_diagrams_heidelberg} are shown together with the models calibrated by each optimization algorithm. The corresponding model parameters, computation times, numbers of objective function evaluations, and relative errors
\begin{equation}
    \mathrm{err}^{z} = \frac{ \left( \sum_{k_0=0}^{M+1-m} \left( q_{\alpha^z,\lambda^z,p^z}^{\,z}(\rho_{k_0}) - q_{k_0}^{z} \right)^2 \right)^{1/2}}{\left( \sum_{k_0=0}^{M+1-m} (q_{k_0}^{z})^2 \right)^{1/2}}, \qquad z\in\{x,y\}
\end{equation}
are summarized in the accompanying tables. 

To compare the performance of all optimization algorithms we evaluate the required running time \texttt{time} and the number of function evaluations $\texttt{feval}$. For CBO and QMC$^2$BO, $\texttt{feval}$ is always given by $\texttt{feval} = N \times (\#\text{Euler-Maruyama steps}+1)$, since we have to fix the number of steps in the beginning.

The complete implementation is available in the GitHub repository accompanying this article \cite{git_2D_data}.

\subsection{Full versus reduced calibration}
\label{sec:numerical_results}

\begin{figure}[htb]
    \centering
    \includegraphics[width=0.9\textwidth]{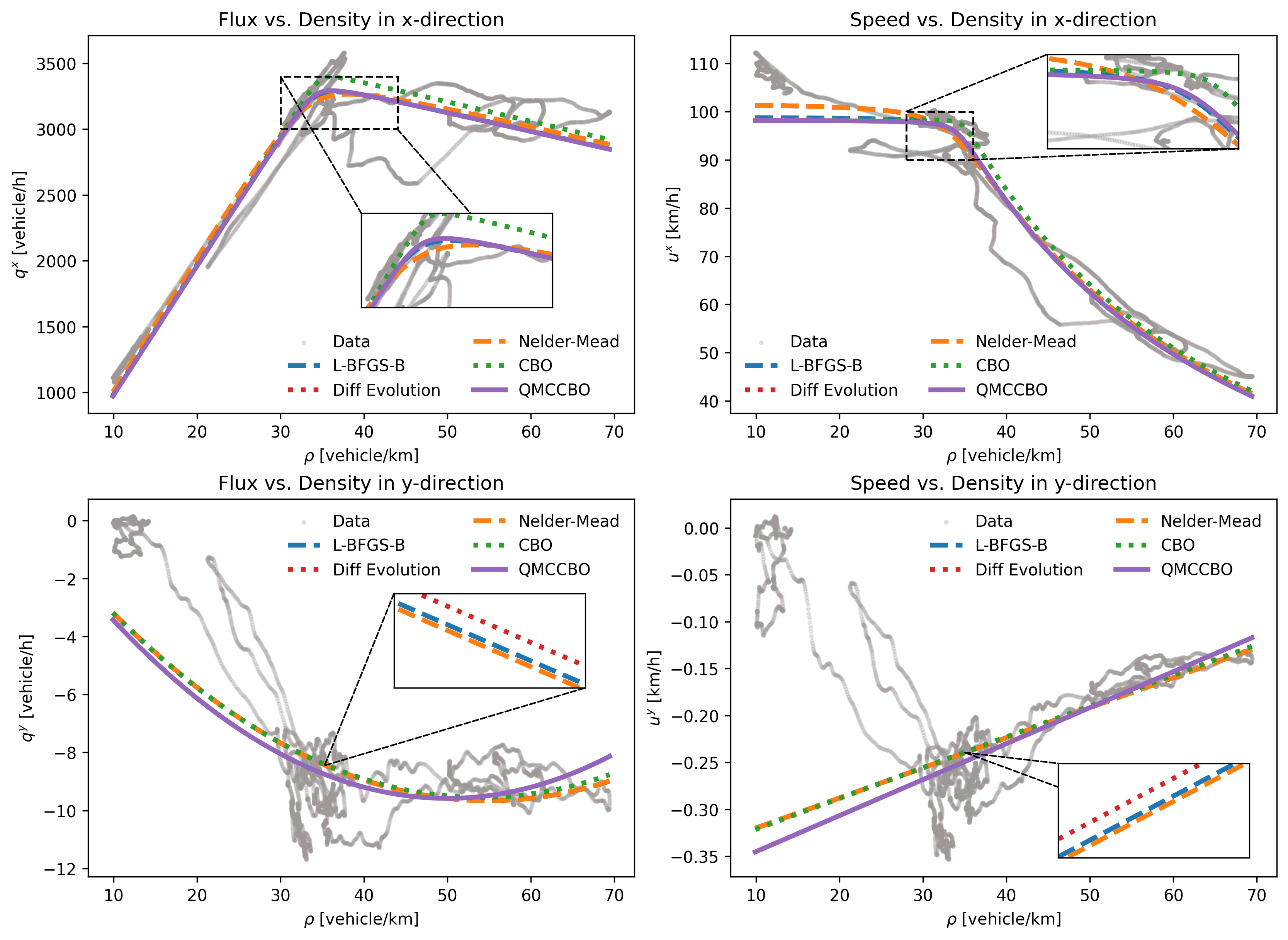}
    \caption{Fundamental (left) and speed-density (right) diagrams obtained from the Off-ramp dataset with full calibration \eqref{eq:full_optimization_variables} solved by each method from Section \ref{sec:optimization}. The upper row shows the movement along the road ($x$-direction), the lower row the movement across lanes ($y$-direction).}
    \label{fig:fundamental_diagrams_offramp_fitted}
\end{figure}

Figure \ref{fig:fundamental_diagrams_offramp_fitted} shows the fundamental diagrams and speed-density diagrams with all fully calibrated models obtained for the Off-ramp dataset. The measured data are shown in light gray and each colored curve corresponds to one of the optimization methods. Since several calibrated models yield similar curves, we zoomed in to selected areas to enable a detailed comparison. The corresponding model parameters, relative errors, computation times, and reported evaluation counters are summarized in Table \ref{tab:parameter_calibration_offramp_full}.

We first consider the movement along the road ($x$-direction). All optimization algorithms recover the structure of the measured data. In particular, they reproduce the linear increase of the flux in the free-flow regime and the flux decreases for the congested system. The models obtained by L-BFGS-B, Differential Evolution and QMC$^2$BO are almost identical. In contrast, CBO predicts the starting of congestion at higher densities than the other models, whereas Nelder-Mead suggests a less sharp turning point.

For the behavior across lanes ($y$-direction), all optimization methods predict a monotone increase of the speed as the traffic density increases, indicating that lane changes towards the right become slower as the road becomes more crowded. The flux, in contrast, is only monotone decreasing up to a turning point between 50 and 60 vehicles per kilometer. This suggests that before the point, more vehicles move towards the rightmost lane as traffic increases. For QMC$^2$BO, this turning point occurs slightly earlier than for the other methods, and afterwards the flux increases faster, indicating that QMC$^2$BO captures a stronger tendency of vehicles to move to the leftmost lane once the road becomes sufficiently congested. This behavior appears plausible, since drivers may want to use all lanes under congestion.

That the model predicted by CBO for the $x$-direction and QMC$^2$BO for the $y$-direction is structurally different is also confirmed through the quantitative results in Table \ref{tab:parameter_calibration_offramp_full}.

\begin{table}[htb]
\centering
\caption{Comparison of the fully calibrated model parameters ($\alpha$, $\lambda$, $p$), relative error ($\texttt{err}$), computation time ($\texttt{time}$) in seconds, and function evaluations ($\texttt{feval}$) for the Off-ramp dataset. The upper table shows the movement along the road ($x$-direction), the lower table the movement across lanes ($y$-direction).}
\small
\begin{tabular}{lcccccc}
\toprule
Method & $\alpha^x$ & $\lambda^x$ & $p^x$ & $\texttt{err}^x$ & $\texttt{time}$ & $\texttt{feval}$\\
\midrule
Nelder-Mead  & 240.49 &  64.65 & 0.13 & $5.593\times10^{-2}$ & 0.51 & 526   \\
L-BFGS-B     & 122.94 & 122.94 & 0.13 & $5.649\times10^{-2}$ & 0.15 & 124   \\
DE           & 240.58 &  64.63 & 0.13 & $5.593\times10^{-2}$ & 3.04 & 3735  \\
CBO          &  52.61 & 287.1  & 0.13 & $6.219\times10^{-2}$ & 2.86 & 99328 \\
QMC$^2$BO    & 100.10 & 150.13 & 0.13 & $5.675\times10^{-2}$ & 2.90 & 99328 \\
\bottomrule
\end{tabular}

\vspace{0.6cm}

\begin{tabular}{lcccccc}
\toprule
Method & $\alpha^y$ & $\lambda^y$ & $p^y$ & $\texttt{err}^y$ & $\texttt{time}$ & $\texttt{feval}$\\
\midrule
Nelder-Mead  &  0.39 & 500.00 & 0.41 & $2.238\times10^{-1}$ & 49.42 & 52367 \\
L-BFGS-B     & 16.71 &  11.78 & 0.42 & $2.239\times10^{-1}$ &  0.06 &    64 \\
DE           &  8.61 &  22.58 & 0.41 & $2.238\times10^{-1}$ &  9.39 & 10080 \\
CBO          &  0.62 & 313.9  & 0.4  & $2.24 \times10^{-1}$  &  3.36 & 99328 \\
QMC$^2$BO    &  0.77 & 284.74 & 0.37 & $2.274\times10^{-1}$ &  3.28 & 99328 \\
\bottomrule
\end{tabular}

\label{tab:parameter_calibration_offramp_full}
\end{table}

In $y$-direction, the relative error of QMC$^2$BO is slightly higher than the other. This suggests that a later turning point is more suitable. In general, the problem in $y$-direction seems more challenging, as all methods have errors of order $10^{-1}$ instead of $10^{-2}$, which is the case for the $x$-direction. Although all optimization algorithms yield similar models for the $x$-direction, CBO has a higher relative error, indicating that the density at which congestion occurs is lower than CBO suggests.

We further observe that the calibrated values of $\alpha^x$ in $x$-direction are strictly positive across all optimization methods. This is precisely what we expected in light of the sign property \eqref{eq:constant_sign} established in Proposition~\ref{prop:flux_structure}.

Interestingly, the calibrated values of $\alpha$ and $\lambda$ vary considerably between methods, whereas the parameter $p$ remains relatively stable at $p^x= 0.13$ and $p^y \in [0.37,0.42]$. Note that the critical density $p^y\rho_{\max}\approx 0.41\cdot 267\approx 109$ vehicles per kilometer lies outside the observed density range of the Off-ramp dataset, so the predicted switch of the lane-changing direction is an extrapolation. Despite the considerable spread in $\alpha$ and $\lambda$, the products $\alpha\lambda$ lie in the much narrower interval $\alpha^x\lambda^x \in [1.502,1.555]\times 10^4$ and $\alpha^y\lambda^y \in [1.944,2.193]\times 10^2$. This is consistent with Proposition \ref{prop:large_lambda_asymptotics} and Corollary \ref{cor:large_lambda_asymptotics_y}, according to which the flux for high $\lambda$ cases depends predominantly on the product $\alpha\lambda$ rather than on $\alpha$ and $\lambda$ individually. In Section \ref{sec:empirical_identifiability}, we will discuss the sensitivity of the full calibration model with respect to these parameters in more detail.

L-BFGS-B uses the fewest function evaluations and has the shortest computation time. Especially for the $y$-direction it is extremely efficient. While the error of Nelder-Mead and Differential Evolution are the lowest, Nelder-Mead is significantly faster and uses fewer function evaluations in $x$-direction. Though, in $y$-direction, the reverse can be observed. This indicates that local optimization can be very efficient when good initialization is chosen. In general, particle based methods need more time and function evaluations.

In general, the computation times of the consensus-based methods should also be interpreted with care. The first three methods rely on highly optimized implementations provided by \texttt{SciPy} \cite{scipy2020}, whereas CBO and QMC$^2$BO are implemented within an external framework.

\medskip

\begin{figure}[ht]
    \centering
    \includegraphics[width=0.9\textwidth]{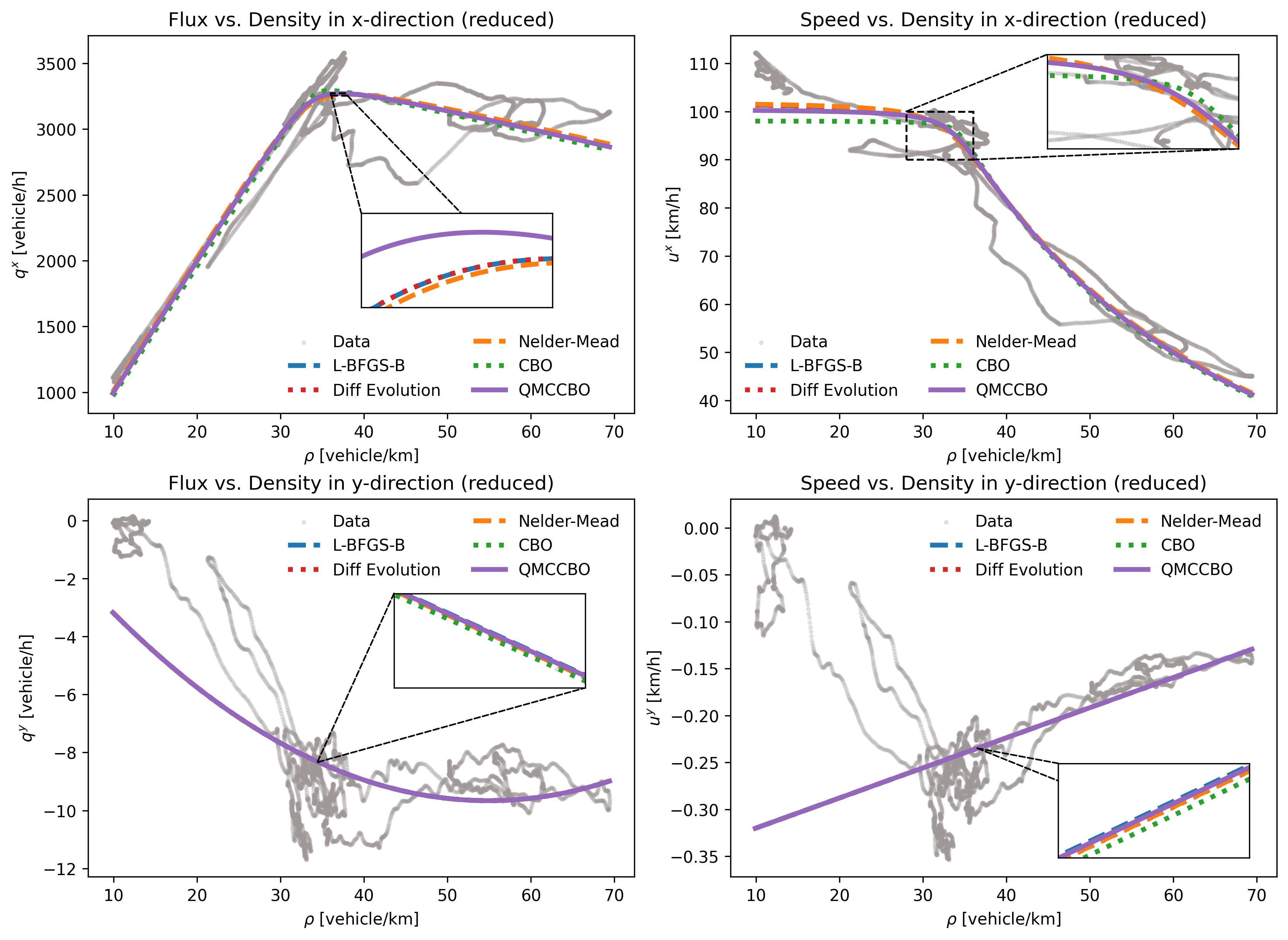}
    \caption{Fundamental (left) and speed-density (right) diagrams obtained from the Off-ramp dataset with reduced calibration \eqref{eq:reduced_optimization_variables} solved by each method from Section \ref{sec:optimization}. The upper row shows the movement along the road ($x$-direction), the lower row illustrates the movement across lanes ($y$-direction).}
    \label{fig:fundamental_diagrams_offramp_fitted_reduced}
\end{figure}

Proposition \ref{prop:alpha_elimination} shows that the reduced problem \eqref{eq:reduced_optimization_variables} is analytically equivalent to the full problem \eqref{eq:full_optimization_variables} with respect to the optimal objective value. Numerically, however, the two formulations may behave very differently. We therefore repeat the experiments for the reduced calibration problem \eqref{eq:reduced_optimization_variables} and compare the result to the full problem. Figure \ref{fig:fundamental_diagrams_offramp_fitted_reduced} shows the corresponding fundamental diagrams and speed-density diagrams.

We observe the same characteristic flux and speed functions as in the full calibration case. However, the curves predicted by the different optimization algorithms now lie closer together than in the full calibration problem. In particular, in the $y$-direction there are hardly any visible differences between the results of the five algorithms.

These observations are consistent with the exact values reported in Table~\ref{tab:parameter_calibration_offramp_reduced}.

\begin{table}[htbp]
\centering
\caption{Comparison of the reduced calibrated model parameters ($\alpha$, $\lambda$, $p$), relative error ($\texttt{err}$), computation time ($\texttt{time}$) in seconds, and function evaluations ($\texttt{feval}$) for the Off-ramp dataset. The upper table shows the movement along the road ($x$-direction), the lower table the movement across lanes ($y$-direction).}
\small
\begin{tabular}{lcccccc}
\toprule
Method & $\alpha^x$ & $\lambda^x$ & $p^x$ & $\texttt{err}^x$ & $\texttt{time}$ & $\texttt{feval}$\\
\midrule
Nelder-Mead  & 248.65 &  62.66 & 0.13 & $5.593\times10^{-2}$ & 0.27 &  215 \\
L-BFGS-B     & 240.80 &  64.57 & 0.13 & $5.593\times10^{-2}$ & 0.08 &   72 \\
DE           & 240.46 &  64.66 & 0.13 & $5.593\times10^{-2}$ & 0.81 &  900 \\
CBO          &  86.8  & 172.8  & 0.13 & $5.691\times10^{-2}$ & 0.19 & 9984 \\
QMC$^2$BO    & 184.99 &  82.99 & 0.13 & $5.606\times10^{-2}$ & 0.20 & 9984 \\
\bottomrule
\end{tabular}

\vspace{0.6cm}

\begin{tabular}{lcccccc}
\toprule
Method & $\alpha^y$ & $\lambda^y$ & $p^y$ & $\texttt{err}^y$ & $\texttt{time}$ & $\texttt{feval}$\\
\midrule
Nelder-Mead  & 0.94 & 205.86 & 0.41 & $2.238\times10^{-1}$ & 0.07 &   64 \\
L-BFGS-B     & 8.21 &  23.69 & 0.41 & $2.238\times10^{-1}$ & 0.08 &   75 \\
DE           & 0.68 & 283.55 & 0.41 & $2.238\times10^{-1}$ & 0.63 &  600 \\
CBO          & 1.09 & 177.88 & 0.41 & $2.238\times10^{-1}$ & 0.13 & 9984 \\
QMC$^2$BO    & 0.46 & 419.81 & 0.41 & $2.238\times10^{-1}$ & 0.17 & 9984 \\
\bottomrule
\end{tabular}

\label{tab:parameter_calibration_offramp_reduced}
\end{table}

In $y$-direction, all methods attain the same error. Although $\alpha^y$ and $\lambda^y$ again differ, which is again consistent with Proposition \ref{prop:large_lambda_asymptotics} and Corollary \ref{cor:large_lambda_asymptotics_y}, all optimization algorithms now agree that $p^y=0.41$ is optimal. Despite using only one tenth of the evaluation budget, the reduced problem reaches comparable accuracy. The same observation holds in the $x$-direction.

Overall, the reduced problem uses less computation time and the number of function evaluations across all optimization algorithms, while still achieving an error that is at least as good as in the full calibration problem. Moreover, it is markedly less sensitive to the choice of optimization algorithm. In $y$-direction, and to a lesser extent in $x$-direction, all methods converge to essentially the same error.

Having established the advantages of the reduced problem, we note that it seems to be easier to solve it. We therefore consider the full calibration problem throughout the remaining analysis, as it is a greater challenge for all five optimization algorithms and therefore allows a more informative comparison between them.

\subsection{Calibration for all datasets}

\begin{figure}[htbp]
    \centering
    \includegraphics[width=0.9\textwidth]{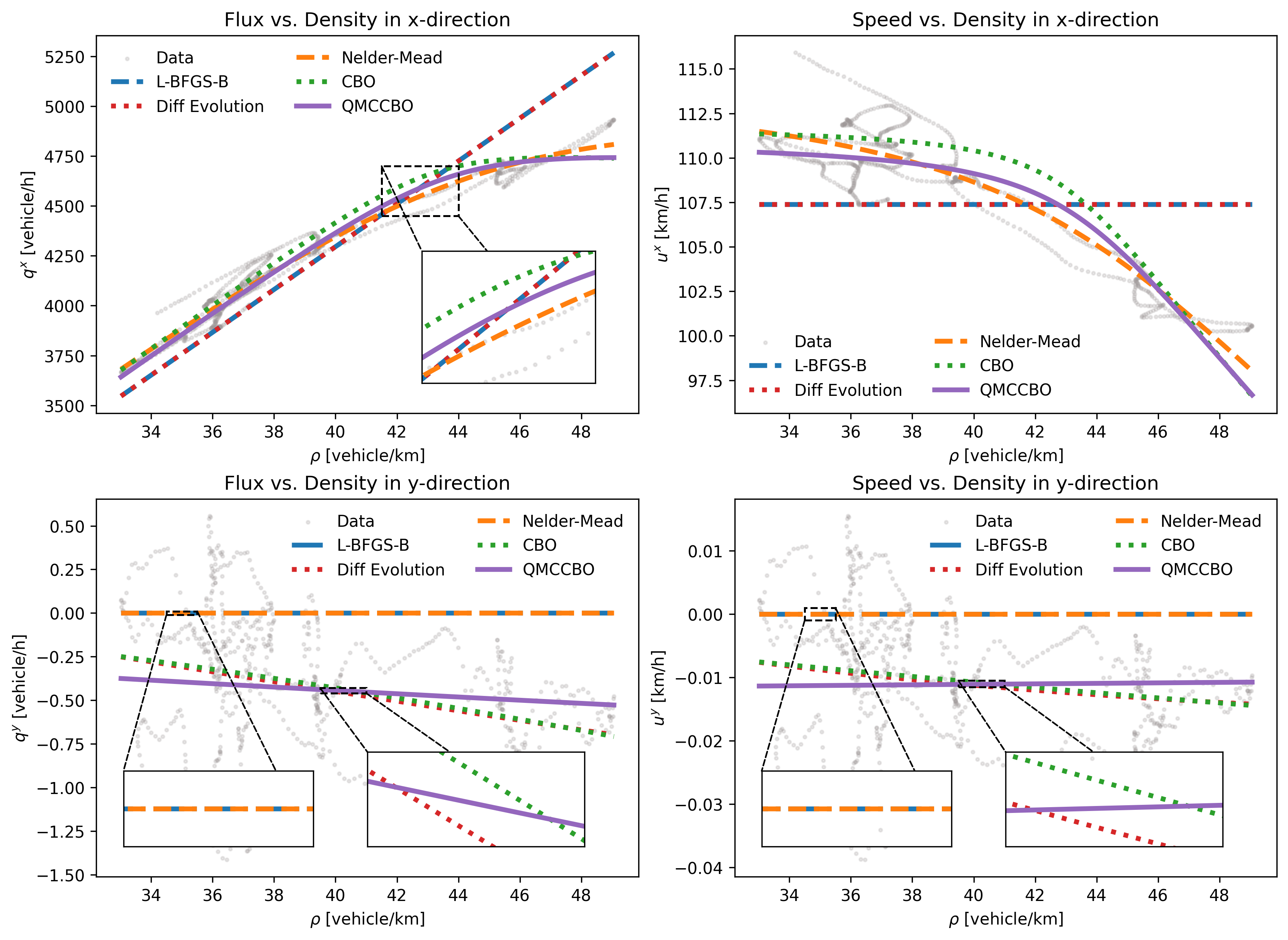}
    \caption{Fundamental (left) and speed-density (right) diagrams obtained from the highD dataset with full parameter calibration problem \eqref{eq:full_optimization_variables} solved by each method from Section \ref{sec:optimization}. The upper row shows the movement along the road ($x$-direction), the lower row the movement across lanes ($y$-direction).}
    \label{fig:fundamental_diagrams_highD_fitted}
\end{figure}

We next consider the full calibration problem of the remaining two datasets. First, we will investigate the first recording of the highD dataset. The fundamental diagrams and speed-density diagrams are shown in Figure \ref{fig:fundamental_diagrams_highD_fitted}. 

The highD dataset contains mostly data-points corresponding to free-flow traffic, which is also reflected in the calibrated models. In $x$-direction, L-BFGS-B and Differential Evolution predict a linear increase of the traffic flux and a constant speed, indicating traffic that is essentially free of congestion. In fact, their relative error is almost three times as big as the one of the other methods. This indicates that the curve shown by Nelder-Mead, CBO, and QMC$^2$BO represents the data better. Here, before a turning point, the system remains in a free-flow regime. After it the flux remains nearly constant and the mean speed decreases linearly with increasing density. This suggests light congestion.

In contrast to the Off-ramp dataset, the highD dataset contains positive and negative flux values in $y$-direction, corresponding to lane changes towards the left and right lane. Nevertheless, all optimization methods predict linear fundamental diagrams for the traffic flux and the speed. In particular, the solutions obtained by L-BFGS-B and Nelder-Mead are constant at zero, indicating no traffic flow at all. Since this clearly contradicts the measured data, these optimization methods fail to capture the observed lane-changing behavior.

The remaining optimization algorithms recover the slight asymmetry behavior in the data. As lane changes towards the rightmost lane occur more frequently, they predict a slightly negative average speed together with a mildly decreasing flux-density relation. Nonetheless, even the best relative error achieved is $\texttt{err}=0.653$, see \cite{git_2D_data}, indicates that the flux model \eqref{eq:flux_param_y} is less suitable for describing the traffic dynamics of the highD dataset, where lane changes in both directions occur frequently.

\medskip

\begin{figure}[htb]
    \centering
    \includegraphics[width=0.9\textwidth]{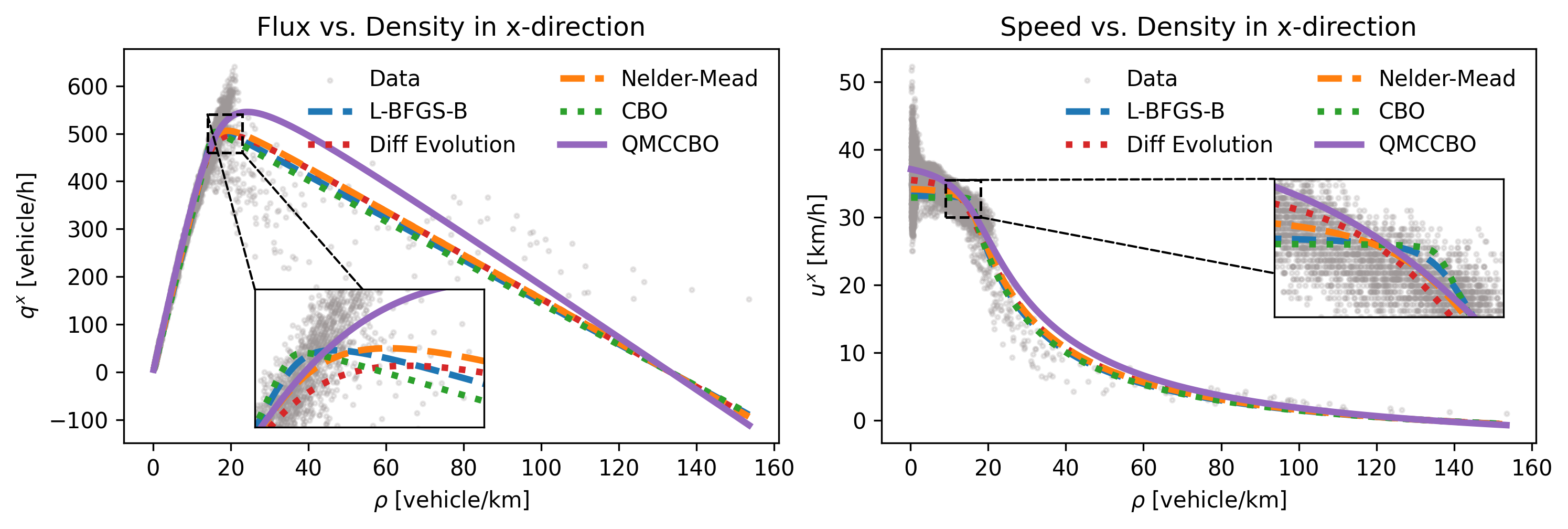}
    \caption{Fundamental (left) and speed-density (right) diagram obtained from the Heidelberg sensor dataset with full parameter calibration problem \eqref{eq:full_optimization_variables} solved by each method from Section \ref{sec:optimization}.}
    \label{fig:fundamental_diagrams_heidelberg_fitted}
\end{figure}

Finally, we consider the Heidelberg sensor dataset. The calibrated fundamental diagrams and speed-density diagrams are shown in Figure \ref{fig:fundamental_diagrams_heidelberg_fitted}.

All optimization methods predict a linear increase of the flux in uncongested situations. The corresponding speed remains constant at a value slightly above the speed limit of $30\,\mathrm{km/h}$. Once congestion starts, the flux decreases linearly and the speed first drops drastically before then slowly approaching its limiting value. Only the solution obtained by QMC$^2$BO exhibits a slightly smoother transition between those two phases, which results in a slightly higher relative error.

All models still underestimate the flux at the highest observed densities. In particular, the fitted flux becomes negative for high densities. This effect is unlikely to be caused by the optimization algorithms themselves. Indeed, the Heidelberg sensor dataset contains mostly observations in the free-flow regime, with comparatively few in highly congested conditions close to $\rho_{\max}$. Consequently, the least-squares objective \eqref{eq:objectives} concentrates more on fitting the densely sampled regions of the fundamental diagram, at the expense of accuracy in the sparsely populated congested regime.

Since the detailed results for the highD and Heidelberg sensor datasets are structurally very similar to Table~\ref{tab:parameter_calibration_offramp_full}, we have already included the relevant aspects when discussing the corresponding figures (Figures~\ref{fig:fundamental_diagrams_heidelberg_fitted} and \ref{fig:fundamental_diagrams_highD_fitted}). Therefore we do not report them explicitly here. Instead, we refer to the GitHub repository accompanying this paper  \cite{git_2D_data}.

For the following three sections, which address the sensitivity and robustness analysis, we choose the Off-ramp dataset as a representative example. Since it provides two-dimensional data with a substantial amount of free-flow and congested observations, it is particularly well suited for examining the behavior of the different optimization algorithms in detail.

\subsection{Sensitivity with respect to optimization variables}
\label{sec:empirical_identifiability}

\begin{figure}[htbp]
    \centering
    \includegraphics[width=1\textwidth]{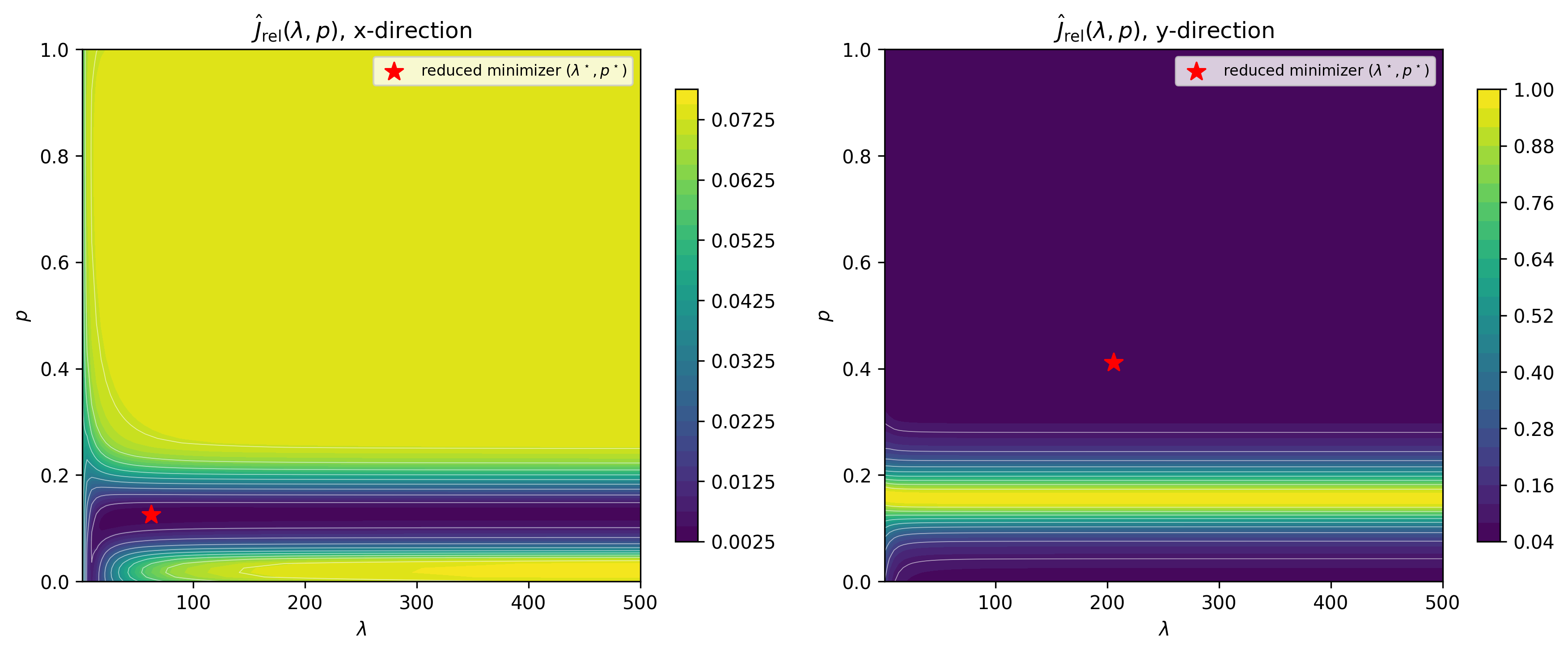}
    \caption{Contour plots of $\widehat J_{\mathrm{rel}}(\lambda,p)$ from \eqref{eq:normalized_reduced_objective} for the Off-ramp dataset. The left panel shows the movement along the road ($x$-direction), the right the movement across lanes ($y$-direction).}
    \label{fig:reduced_objective_contours_offramp}
\end{figure}

Consistent with Proposition~\ref{prop:large_lambda_asymptotics}, the numerical experiments show that the flux is mostly determined by $p$ and the product $\alpha\lambda$ rather than by $\alpha$ and $\lambda$ individually. Since the reduced calibration problem \eqref{eq:reduced_calibration} is equivalent to the full problem \eqref{eq:objectives}, and since the optimal $\alpha$ can be uniquely recovered from $\lambda$ and $p$ via \eqref{eq:alpha_star}, it suffices to investigate the sensitivity of the calibration problem with respect to $\lambda$ and $p$ alone. We therefore take a closer look at the sensitivity of the full calibration problem with respect to these two parameters by considering the normalized reduced objective
\begin{equation}
\label{eq:normalized_reduced_objective}
    \widehat J_{\mathrm{rel}}(\lambda,p) :=
    \frac{\widehat J(\lambda,p)}{\frac12\sum_{k_0=0}^{M+1-m} (q^z_{k_0})^2},
\end{equation}
where $\widehat J$ is the reduced least-squares functional from Proposition~\ref{prop:alpha_elimination}. Figure \ref{fig:reduced_objective_contours_offramp} shows the resulting contour plots for the Off-ramp dataset.

In both panels, the objective $\widehat J_{\mathrm{rel}}(\lambda,p)$ is flat for different $\lambda$ around the optimal value of $p$. In addition, slight changes of $p$, while remaining $\lambda$ fixed, may increase the function value drastically. This confirms that $p$ is well identified by the data, whereas $\lambda$ is flexible. This is due to the possibility to adjust the value of $\alpha\lambda$ by varying $\alpha$. 

In $x$-direction, we observe that low objective values are attained only in a very narrow range of $p^x$ values, namely $p^x\approx 0.13$, and that the function increases sharply as soon as we move away from this range. In $y$-direction, the optimal $p^y$ value lies in a considerably flatter region. This supports our earlier observations, where all optimization methods agreed on $p^x=0.13$, while in the $y$-direction, $p^y \in [0.37,0.42]$.

We now quantify this behavior more precisely. At an interior minimizer $(\lambda^\star,p^\star)$, a small eigenvalue of the Hessian
\begin{equation}
\label{eq:reduced_hessian}
    H_{\mathrm{red}} := \nabla^2_{(\lambda,p)}\widehat J(\lambda^\star,p^\star)
\end{equation}
means that the objective function is relatively flat in the direction of the corresponding eigenvector, while a large one indicates a steep increase. We evaluate $H_{\mathrm{red}}$ at the minimizer $(\lambda^\star,p^\star)$ of the reduced problem obtained by Nelder-Mead, since this method attains the lowest error in both directions. This minimizer lies in the interior of the feasible set for both directions. Following the central-difference construction outlined by \citenote{remark following eq. (8.20)}{nocedal2006numerical}, we approximate the entries of $H_{\mathrm{red}}$ using symmetric perturbations with step size $h=10^{-3}$ in both variables.

Since $\lambda\in[0,500]$ and $p\in[0,1]$ live on very different scales, the eigenvalues of $H_{\mathrm{red}}$ and their ratio depend strongly on the parametrization. We therefore consider the scaled variable $\tilde\lambda=\lambda/500$, for which the Hessian becomes $\widetilde H_{\mathrm{red}}=D\,H_{\mathrm{red}}\,D$ with $D=\mathrm{diag}(500,1)$. The resulting eigenvalues $\mu_1\le\mu_2$ and condition numbers are summarized in Table~\ref{tab:hessian_eigenvalues}.

\begin{table}[h]
\centering
\caption{Eigenvalues $\mu_1\le\mu_2$ of the reduced Hessian $H_{\mathrm{red}}$ from \eqref{eq:reduced_hessian} in the scaled variables $(\tilde\lambda,p)$, evaluated at the Nelder-Mead minimizer of the reduced problem for the Off-ramp dataset, and their ratio $\frac{\mu_2}{\mu_1}$.}
\begin{tabular}{lccc}
\toprule
Direction & $\mu_1$ & $\mu_2$ & $\frac{\mu_2}{\mu_1}$ \\
\midrule
$x$ & $1.70\times10^{9}$ & $6.66\times10^{11}$ & $3.93\times10^{2}$ \\
$y$ & $4.36\times10^{1}$ & $4.62\times10^{5}$  & $1.06\times10^{4}$ \\
\bottomrule
\end{tabular}

\label{tab:hessian_eigenvalues}
\end{table}

In $x$-direction, the eigenvalue $\mu_2$, associated with $p^x$, is extremely large, reflecting the sharp valley around $p^x\approx0.13$ observed in the contour plot. Even after scaling, the condition number is of order $4\times10^{2}$, so the objective is still considerably flatter along $\lambda^x$ than along $p^x$.

In $y$-direction, $\mu_2$ is much smaller than in $x$-direction, which is consistent with the flatter region around the optimal $p^y$. At the same time, $\mu_1$, associated with $\lambda^y$, is much smaller as well, so the flatness in $\lambda$-direction is even more pronounced. This results in a condition number of $1.06\times10^{4}$. The poorer conditioning in $y$-direction is therefore mainly due to the very flat $\lambda$-direction rather than to the $p$-direction.

\subsection{Robustness of the optimization results}
\label{sec:optimization_robustness}

To assess the robustness of the different optimization algorithms, we perform a Monte-Carlo analysis on the full calibration problem \eqref{eq:full_optimization_variables} consisting of $1000$ independent runs. Throughout this section, all algorithmic parameters are kept identical to those used in Section \ref{sec:numerical_results}.

Since Nelder-Mead, L-BFGS-B, and QMC$^2$BO are deterministic, they always return the same solution. Consequently, the Monte Carlo analysis is primarily concerned with the stochastic optimization methods, Differential Evolution and CBO. For each method, Table \ref{tab:monte_carlo_errors} reports the average relative error and its empirical standard deviation, denoted by $\texttt{dev}$.

\begin{table}[ht]
\centering
\caption{Comparison of the mean relative error ($\texttt{err}$) and its standard deviation ($\texttt{dev}$) obtained from $1000$ Monte Carlo simulations on the Off-ramp dataset for each optimization method from Section \ref{sec:optimization}. The results are reported for both the movement along the road ($x$-direction) and across lanes ($y$-direction).}
\setlength{\tabcolsep}{12pt}

\begin{tabular}{lcccc}
\toprule
Method 
& $\texttt{err}^x$ 
& $\texttt{dev}^x$ 
& $\texttt{err}^y$ 
& $\texttt{dev}^y$ \\
\midrule

Nelder-Mead 
& $5.593\times10^{-2}$ 
& $0.00$ 
& $2.238\times10^{-1}$ 
& $0.00$ \\

L-BFGS-B
& $5.649\times10^{-2}$ 
& $0.00$ 
& $2.239\times10^{-1}$ 
& $0.00$ \\

DE
& $5.928\times10^{-2}$ 
& $2.64\times10^{-2}$ 
& $2.264\times10^{-1}$ 
& $9.88\times10^{-3}$ \\

CBO
& $8.066\times10^{-2}$ 
& $1.35\times10^{-1}$ 
& $5.809\times10^{-1}$ 
& $9.3\times10^{-1}$ \\

QMC$^2$BO
& $5.675\times10^{-2}$ 
& $0.00$ 
& $2.274\times10^{-1}$ 
& $0.00$ \\

\bottomrule
\end{tabular}

\label{tab:monte_carlo_errors}

\end{table}

We observe that Differential Evolution is quite robust, while CBO shows deviations in between runs. Although the average error remains of the same order, $10^{-2}$ in the $x$-direction and $10^{-1}$ in the $y$-direction, as observed in Section~\ref{sec:numerical_results}, Differential Evolution is particularly robust, with a deviation of order $10^{-2}$ in the $x$-direction and $10^{-3}$ in the $y$-direction, indicating that almost all runs converge to very similar values and that the relative error remains stable across repetitions.

\subsection{Dependence on the initial configuration}
\label{sec:initialization_sensitivity}

The Monte Carlo study in Section~\ref{sec:optimization_robustness} showed that stochastic optimization methods may vary between runs. We now investigate whether a similar sensitivity is also present for deterministic optimization methods with respect to their initialization. For each optimization algorithm, three different initial configurations are considered and the resulting relative errors for the Off-ramp dataset are summarized in Table \ref{tab:initialization_sensitivity}.

Since all optimization methods use different optimization strategies, the initial configuration differ structurally from algorithm to algorithm. For Nelder-Mead, the initialization consists of an initial simplex, i.e., a set of $D+1=4$ vertices in $\mathbb{R}^3$. We consider the three simplices
$$
\begin{aligned}
s_1=\{&
(-100,\,0,\,0.9),\,
(0,\,100,\,0.9),\,
(10,\,0,\,0.9),\,
(-400,\,0,\,0.6)
\},\\[1ex]
s_2=\{&
(0,\,250,\,0.5),\,
(50,\,150,\,0.5),\,
(0,\,350,\,0.5),\,
(0,\,250,\,0.9)
\},\\[1ex]
s_3=\{&
(-450,\,10,\,0.05),\,
(450,\,200,\,0.15),\,
(-150,\,340,\,0.35),\,
(-450,\,10,\,0.15)
\}.
\end{aligned}
$$
For L-BFGS-B, we use the three initial points
$$
p_1=(100,100,0.5),\qquad p_2=(10,5,0.5),\qquad p_3=(-30,150,0.7).
$$
For Differential Evolution and CBO, three independent runs with different uniformly distributed initial configurations are conducted. Finally, for QMC$^2$BO, we replace the Sobol \cite{Sobol_points} initialization used in Section \ref{sec:numerical_results} by two different quasi-Monte Carlo point sets, namely Halton points \cite{Halton_points} and rank-$1$ lattice point sets \cite{Lattice_points}.

\begin{table}[ht]
\centering
\caption{Comparison of the relative errors ($\texttt{err}$) obtained from the Off-ramp dataset for different initialization strategies of each method from Section~\ref{sec:optimization}. For Nelder-Mead and L-BFGS-B, different initial simplices and initial values are considered, for Differential Evolution and CBO independent random runs and for QMC$^2$BO different quasi-Monte Carlo sampling strategies.}
\setlength{\tabcolsep}{10pt}

\begin{tabular}{llcc}
\toprule
Method & Initialization & $\texttt{err}^x$ & $\texttt{err}^y$\\
\midrule

\multirow{3}{*}{Nelder-Mead}
& $s_1$ & $2.708\times10^{-1}$ & $1.0\times10^{0}$\\
& $s_2$ & $5.593\times10^{-2}$ & $2.238\times10^{-1}$\\
& $s_3$ & $1.0\times10^{0}$    & $1.0\times10^{0}$\\

\midrule

\multirow{3}{*}{L-BFGS-B}
& $p_1$   & $5.649\times10^{-2}$ & $1.0\times10^{0}$\\
& $p_2$   & $2.745\times10^{-1}$  & $2.239\times10^{-1}$\\
& $p_3$   & $1.0\times10^{0}$    & $1.0\times10^{0}$\\

\midrule

\multirow{3}{*}{DE}
& Run 1 & $5.632\times10^{-2}$ & $2.649\times10^{-1}$\\
& Run 2 & $5.595\times10^{-2}$ & $2.239\times10^{-1}$\\
& Run 3 & $5.601\times10^{-2}$ & $2.256\times10^{-1}$\\

\midrule

\multirow{3}{*}{CBO}
& Run 1 & $5.64\times10^{-2}$ &  $2.627\times10^{-1}$\\
& Run 2 & $5.663\times10^{-2}$ & $1.0\times10^{0}$\\
& Run 3 & $5.649\times10^{-2}$ & $2.242\times10^{-1}$\\

\midrule

\multirow{3}{*}{QMC$^2$BO}
& Sobol      & $5.675\times10^{-2}$ & $2.274\times10^{-1}$\\
& Halton     & $5.647\times10^{-2}$ & $2.291\times10^{-1}$\\
& Lattice    & $5.855\times10^{-2}$ & $2.318\times10^{-1}$\\

\bottomrule
\end{tabular}

\label{tab:initialization_sensitivity}

\end{table}
These results confirm the Monte-Carlo analysis from the previous section. Differential Evolution is stable in both directions. CBO exhibits for the $y$-direction one run that does not appear to converge at all. It is precisely such runs that explain the high deviation of CBO compared to Differential Evolution.

We also observe that local optimization methods are sensitive to initialization. Nelder-Mead and L-BFGS-B each have one initialization for which the estimate in the $x$-direction is close to the optimum, one that is further off by an order of $10^{-1}$, and one that does not converge at all, resulting in an error of $1$. For the $y$-direction, the situation is even worse. Here, for each method, two of the initial values fail to converge at all.

The failures with $\texttt{err}=1$ have a structural cause. Since $g_{0,p}\equiv0$, the model flux vanishes for $\lambda=0$ and every $\alpha,p$, so the relative error equals one. As $g_{\lambda,p}$ depends on $\lambda$ only through $\lambda^2$, the plane $\{\lambda=0\}$ consists of stationary points of $J^z$. Local methods that reach this plane, or start on it, terminate there without indication of failure. For example, all vertices of the simplex $s_1$ satisfy $\lambda=0$ or $\alpha=0$. The failed runs are therefore a consequence of this plateau and unfavorable starting points rather than a general deficiency of the methods.
In contrast, QMC$^2$BO shows a robust behavior with respect to initialization. All three quasi-Monte Carlo point sets produce similar errors, indicating that the quality of convergence is mostly determined by the uniform coverage rather than by the specific construction of the point set.

\section{Conclusion}
\label{sec:conclusion}

In this work, we calibrated a two-dimensional first-order traffic flow model using three real-world datasets. After transforming the observations into comparable macroscopic fundamental diagrams, we analyzed the parameterized flux family and the associated least-squares problem. We then solved the calibration problems with five different optimization algorithms.

We observed in Proposition~\ref{prop:wellposedness} that the underlying conservation law is well-posed in the Kružkov sense. For the parameterized flux, Proposition~\ref{prop:flux_structure} shows that the sign of the curvature parameter $\lambda$ cannot be identified and that every non-degenerate flux is strictly concave or convex, with a sign determined by $\alpha$. Consistently with this, all optimization methods returned $\alpha^x>0$ for the movement along the road. The constant-sign property motivates to consider a different flux model for the movement across lanes, namely \eqref{eq:flux_param_y}.

In Proposition~\ref{prop:alpha_elimination} and Corollary~\ref{cor:alpha_elimination_y}, we formulate a simplified calibration problem, which reduces the three-dimensional problem exactly to a two-dimensional problem in $(\lambda,p)$. Although the reduced problem is analytically equivalent, it performs considerably better in the numerical experiments. All methods required fewer function evaluations and less computation time to reach the same or even a smaller error. For instance, Nelder-Mead required $64$ instead of $52367$ evaluations in the $y$-direction for the Off-ramp dataset. Moreover, the predictions became nearly independent of the optimization method, which indicated that this problem is more robust.

Furthermore, in Proposition~\ref{prop:large_lambda_asymptotics} and Corollary~\ref{cor:large_lambda_asymptotics_y} we established that for large $\lambda$ the flux depends essentially only on $p$ and on the product $\alpha\lambda$. Both effects are also visible in the experiments, where $p$ and $\alpha\lambda$ remain within narrow ranges, while $\alpha$ and $\lambda$ itself may vary greatly. Additionally, the sensitivity analysis in Section~\ref{sec:empirical_identifiability} quantifies this behavior. The reduced objective is almost flat along $\lambda$ but increases steeply in $p$.

Apart from individual exceptions, such as the essentially linear fluxes returned by L-BFGS-B and Differential Evolution for the highD dataset, all optimization methods yielded consistent and accurate fits. The only systematic difference concerns the driving direction. In $x$-direction, the relative errors are of order $10^{-2}$, in $y$-direction they are of order $10^{-1}$. Since this difference occurs uniformly across all optimization methods, this is most likely caused by properties of the data rather than of the optimization algorithms. This motivates future work on more flexible flux models for the movement across the road.

Regarding the algorithms, local methods were the most efficient when properly initialized. L-BFGS-B required the fewest function evaluations, and Nelder-Mead, together with Differential Evolution, attained the smallest errors. However, both local methods strongly depend on their initialization and failed for several initial configurations, caused by the stationary plateau $\{\lambda=0\}$ on which the model flux vanishes. Differential Evolution was more reliable than CBO, which achieved comparable errors in most runs but occasionally failed to converge. QMC$^2$BO is deterministic and therefore avoids this variability. It yielded similar errors for different Quasi-Monte Carlo initialization.

\section*{Use of Generative-AI tools declaration}
The authors declare they have used Artificial Intelligence (AI) tools in the creation of
this article.

\medskip\noindent
\textbf{AI tools used:} ChatGPT (OpenAI).

\medskip\noindent
\textbf{How were the AI tools used?} ChatGPT was used for language editing and for
reviewing parts of the mathematical exposition, in particular to identify potential
inconsistencies and to suggest clearer formulations. It was also used to assist with
parts of the Python implementation of the data preprocessing and the numerical
experiments. No AI tool was used to generate or modify the traffic data, and all
numerical results were obtained with code that was reviewed and tested by the authors.
The authors are aware that generative AI tools may produce inaccurate or biased output;
all AI-generated suggestions and code were therefore checked by the authors, and all
references were verified against the original sources. The authors take full
responsibility for the content of this article.

\medskip\noindent
\textbf{Where in the article is the information located?} Language editing: entire
manuscript. Review of the mathematical exposition: Section~\ref{sec:model}. Code
assistance: implementation underlying Sections~\ref{sec:data_preparation}
and~\ref{sec:numerical_experiments}, available in the repository~\cite{git_2D_data}

\section*{Acknowledgments}
The authors gratefully acknowledge Andreas Neuenkirch for insightful discussions and valuable suggestions that helped improve this work.

\section*{Conflict of interest}
All authors declare no conflicts of interest in this paper.

\bibliographystyle{siam}
\bibliography{literature_revised_v17}

\end{document}